%% file: main.tex
\documentclass[11pt,a4paper]{article}
\usepackage{cmap} 
\usepackage[utf8]{inputenc}				
\usepackage[T1]{fontenc}                

\usepackage{a4wide}
\usepackage{amssymb,amsmath,amsthm,bm}
\usepackage{mathtools}
\usepackage{xifthen}
\usepackage[english]{babel}
\usepackage{epsfig}
\usepackage{tikz}
\usetikzlibrary{arrows}
\usepackage{subfigure}
\usepackage{enumitem} 
\usepackage{graphicx}
\usepackage{todonotes}
\usepackage{standalone}
\usepackage{fullpage}

\usepackage{hyperref}

\newtheorem{theorem}{Theorem}

\newtheorem{lemma}[theorem]{Lemma}

\newtheorem{claim}[theorem]{Claim}

\newtheorem{construction}{Construction}
\newtheorem{problem}[theorem]{Problem}
\newtheorem{obs}[theorem]{Observation}
\newtheorem{defi}[theorem]{Definition}
\newtheorem{definition}[theorem]{Definition}

\usepackage{indentfirst}
\usepackage{xspace} 

\def\R{\mbox{\ensuremath{\mathbb R}}\xspace}

\def\C{\mbox{\ensuremath{\mathcal C}}\xspace}
\def\F{\mbox{\ensuremath{\mathcal F}}\xspace}
\def\G{\mbox{\ensuremath{\mathcal G}}\xspace}

\def\F{\mbox{\ensuremath{\mathcal F}}\xspace}
\def\P{\mbox{\ensuremath{\mathcal P}}\xspace}
\def\S{\mbox{\ensuremath{\mathcal S}}\xspace}

\newcommand{\pointSet}{P}
\newcommand{\point}[1][]{p_{#1}}
\newcommand{\coords}[1][]{(x_{#1},y_{#1})}
\newcommand{\coordsP}[1][]{(x'_{#1},y'_{#1})}
\newcommand{\region}[1][]{R_{#1}}

\newcommand{\lo}[1]{\underline{#1}}
\newcommand{\hi}[1]{\overline{#1}}

\DeclareMathOperator{\image}{image}

\DeclareMathOperator{\sat}{sat}
\DeclareMathOperator{\osat}{ssat}
\DeclareMathOperator{\ram}{ram}
\DeclareMathOperator{\conv}{conv}
\DeclareMathOperator{\interior}{int}

\newcommand{\lp}{\left (}
\newcommand{\rp}{\right )}

\newcommand{\abs}[1]{\left\lvert{#1}\right\rvert}
\newcommand{\floor}[1]{\left\lfloor{#1}\right\rfloor}
\newcommand{\ceil}[1]{\left\lceil{#1}\right\rceil}

\begin{document}

	\title{Saturation problems in the Ramsey theory of graphs, posets and point sets}

	\author{Gábor Damásdi\thanks{MTA-ELTE Lend\"ulet Combinatorial Geometry Research Group, Budapest. The project was  supported by the European Union, co-financed by the European Social Fund (EFOP-3.6.3-VEKOP-16-2017-00002).
	} \and
	Balázs Keszegh\thanks{Alfréd Rényi Institute of Mathematics and MTA-ELTE Lend\"ulet Combinatorial Geometry Research Group, Budapest. Research supported by the by the Lend\"ulet program of the Hungarian Academy of Sciences (MTA), under the grant LP2017-19/2017 and by the National Research, Development and Innovation Office -- NKFIH under the grant K 116769 and K 132696.} \and
	David Malec\thanks{Department of Economics, University of Maryland, College Park, Maryland.} \and
	Casey Tompkins\thanks{Discrete Mathematics Group, Institute for Basic Science (IBS), Daejeon, Republic of Korea. This author was supported by the grant IBS-R029-C1.} \and
	Zhiyu Wang\thanks{University of South Carolina, Columbia, SC, 29208, (zhiyuw@math.sc.edu). This author was supported in part by NSF grant DMS-1600811.}\and
	Oscar Zamora\thanks{Central European University, Budapest. Universidad de Costa Rica, San Jos\'e. The research of this author was partially supported by the National Research, Development and Innovation Office NKFIH, grant  K132696.}}
	
	\maketitle
	
\begin{abstract}
In 1964, Erd\H{o}s, Hajnal and Moon introduced a saturation version of Tur\'an's classical theorem in extremal graph theory. In particular, they determined the minimum number of edges in a $K_r$-free, $n$-vertex graph with the property that the addition of any further edge yields a copy of $K_r$. We consider analogues of this problem in other settings. We prove a saturation version of the Erd\H{o}s-Szekeres theorem about monotone subsequences and saturation versions of some Ramsey-type theorems on graphs and Dilworth-type theorems on posets.  

We also consider \begin{em}semisaturation\end{em}  problems, wherein we allow the family to have the forbidden configuration, but insist that any addition to the family yields a new copy of the forbidden configuration.  In this setting, we prove a semisaturation version of the Erd\H{o}s-Szekeres theorem on convex $k$-gons, as well as multiple semisaturation theorems for sequences and posets.
\end{abstract}


\section{Introduction}
Extremal problems have a long history in combinatorics originating with the results of Mantel~\cite{Mantel} in 1907 and Turán~\cite{Turan} in 1947 determining the maximum number of edges in a triangle- and $K_r$-free, $n$-vertex graph, respectively. Erdős, Hajnal and Moon~\cite{ehm} investigated the dual problem, called the saturation problem, wherein one aims to minimize the number of edges in a $K_r$-free, $n$-vertex graph, such that the addition of any edge yields a copy of $K_r$. Since their initial result, the saturation problem has been considered for a variety of graphs.  Of particular note is a theorem of Kászonyi and Tuza~\cite{Tuza}, which showed that for any graph $H$, the minimum number of edges in an $H$-saturated, $n$-vertex graph is at most linear in $n$.  

Going beyond graphs, saturation problems have been considered in several other settings.  A structure which is maximal with respect to some property is said to \emph{saturate} that property. A maximum size saturating structure is called an \emph{extremal structure}, while a minimum size saturating structure is called a \emph{minimal saturating structure}. For intersecting hypergraphs, a saturation version of the Erdős-Ko-Rado theorem~\cite{ekr} (uniform setting) was proven by Füredi~\cite{furedi}.  In particular, he showed that for a given uniformity $r$, there exists a family of approximately $3r^2/4$ sets of size $r$ with the property that adding any further set yields a pair of disjoint sets, disproving a conjecture of Meyer~\cite{meyer}. In the nonuniform setting, it is well known that all maximal intersecting families of subsets of an $n$ element set have the same size, namely $2^{n-1}$. This result was extended to the case of families without $k$-matchings by Buci\'c~\emph{et al.}~\cite{bucic}.
In the setting of forbidden (induced or non-induced) posets in the Boolean lattice the saturation problem has been investigated by Ferrara~\emph{et al.}~\cite{Ferrara}, and further results in this direction were obtained in~\cite{mar} and~\cite{satsolved}.

Parallel with the development of extremal combinatorics, Ramsey theory has been investigated extensively.  This topic begins with the seminal result of Ramsey~\cite{ramsey}, which states that for any integers $c,r,k$ there is an integer $N$ such that any $c$-coloring of the edges of an $r$-uniform hypergraph on $N$ vertices contains a monochromatic complete graph of size $k$ in some color.  This initial result gave rise to a variety of problems where in place of a complete hypergraph one is given hypergraphs $F_1,F_2,\dots,F_c$, and seeks to minimize the value of $N$ which yields, for all $c$-colorings of the complete $r$-uniform, $N$-vertex hypergraph, a copy of $F_i$ in color $i$ for some $i=1,2,\dots,c$.

Ramsey-type problems may be interpreted as extremal problems in the following way. One wishes to maximize the number of vertices $n$ in such a way that there exists a coloring, such that for all $i$, we find no copy of $F_i$ in color $i$ (so $n=N-1$, where $N$ is defined as above). With this interpretation, it becomes natural to ask the corresponding minimal saturation problem, wherein we seek to minimize $n$ such that the hypergraph can be $c$-colored without a monochromatic copy of $F_i$ in color $i$, but if we extend this $c$-colored hypergraph to a $c$-colored hypergraph on $n+1$ vertices, then we have for some $i$, a monochromatic copy of some $F_i$ in color $i$. 

Finally, we mention that many classical results, such as Dilworth's theorem on posets and the Erd\H{o}s-Szekeres theorem for sequences or cups and caps, can be interpreted as Ramsey-type problems where the allowed colorings of the hypergraph are restricted in some way.  As such, we may again consider the corresponding saturation versions of these results. In this paper, we initiate such a study of Ramsey-type saturation problems. We concentrate on well-known settings (graphs, posets, monotone and convex subsets of point sets), and in several cases we manage to prove tight bounds.

In addition, we also consider the corresponding semisaturation problems, a notion introduced by F\"{u}redi and Kim~\cite{semi} (also sometimes called oversaturation or strong saturation). In the graph setting the semisaturation problem is the following:  Given a graph $F$, what is the minimum number of edges in an $n$-vertex graph $G$ with the property that adding any edge to $G$ yields a copy of $F$ containing that edge.  Note that now we allow the graph $G$ to contain $F$ as a subgraph. We will consider semisaturation problems for sequences, cups and caps, posets and point sets as well as for the Ramsey problem on graphs.  Note that by definition, the semisaturation number is always at most the saturation number which is in turn at most the extremal (Ramsey) number.

In the rest of this section we provide a precise formulation of each of the saturation and semisaturation problems that are considered in the paper and our results for each case. We also briefly summarize the known results about the corresponding extremal problem in order to contrast them with our minimal saturation results. Sections~\ref{section:graph}--\ref{section:conv} contain the proofs of these results. Finally, in Section~\ref{general} we rigorously define the general framework that was hinted at above and illustrate how these problems fit into it. 

\subsection*{Graphs}

 Let $\G$ be the family of (labeled) complete graphs whose edges are colored with $c$ colors (numbered by $1,2,\dots, c$). Given $G, G'\in \G$, we say $G'$ \emph{extends} $G$ if $G$ is a proper subgraph of $G'$, i.e., $G'$ can be obtained from $G$ by iteratively adding a new vertex and colored edges connecting the new vertex with each of the existing vertices. 
 A member $G$ of $\G$ is called \emph{$(k_1, k_2, \ldots, k_c)$-saturated} if for every $i\in [c]$, the graph $G$ does not contain a monochromatic $K_{k_i}$ of color $i$, but every $G' \in \G$ that extends $G$ contains a monochromatic $K_{k_i}$ of color $i$ for some $i$. A graph $G \in \G$ is called \emph{$(k_1, k_2, \ldots, k_c)$-semisaturated} if for every $G'\in \G$ that extends $G$, there exists some $i\in [c]$ such that $G'$ contains a copy of a monochromatic $K_{k_i}$ of color $i$ which is not in $G$.

Clearly the size of the largest $(k_1, k_2, \ldots, k_c)$-saturated graph in $\G$, which we denote by $\ram_{\G}(k_1,\dots, k_c)$, is equal to the usual Ramsey number minus one. Let $\sat_{\G}(k_1,\dots k_c)$ denote the size of the smallest saturated $G\in \G$, and finally let $\osat_{\G}(k_1,\dots k_c)$ denote the size of the smallest $(k_1, \dots, k_c)$-semisaturated $G\in \G$. 
From the definition it is clear that \[\osat_{\G}(k_1,\dots, k_c)\le \sat_{\G}(k_1,\dots, k_c) \le \ram_{\G}(k_1,\dots, k_c).\] 
For convenience, we also use $\sat_{\G}(k;c)$ to denote $\sat_{\G}(k_1, k_2, \dots, k_c)$ and $\osat_{\G}(k;c)$ to denote $\osat_{\G}(k_1, k_2, \dots, k_c)$ when $k_1 = k_2 = \cdots = k_c = k$.

For a fixed $k$ and growing $l$, the following results about Ramsey numbers are known (due to Bohman and Keevash~\cite{BK} and Ajtai, Komlós and Szemerédi~\cite{AKS}, respectively): 
\[c'_k\frac{l^{\frac{k+1}{2}}}{(\log l)^{\frac{k+1}{2}-\frac{1}{k-2}}}\le \ram_{\G}(k,l)\le c_k\frac{l^{k-1}}{(\log l)^{k-2}}.\]
For the case $l=3$, it is known that 
\[\ram_{\G}(k,3)=\Theta\left(\frac{k^2}{\log k}\right).\] 
The upper bound was proven by Ajtai, Komlós, and Szemerédi~\cite{AKS}; the lower bound was obtained originally by Kim~\cite{Kim}, and subsequently improved by Fiz Pontiveros, Griffiths and Morris~\cite{FGM} and Bohman and Keevash~\cite{BK}.



We prove the following results:

\begin{theorem}\label{theorem:satgraphs}
    For two colors,
	 \[\osat_{\G}(k,l)= \sat_{\G}(k,l)= (k-1)(l-1),\]
	and for more than two colors, 
	\[(k_1-1)(k_2+\dots+k_c-2c+3)\le \osat_{\G}(k_1, \ldots, k_c)\le \sat_{\G}(k_1, \ldots, k_c)\le (k_1-1)\cdots (k_c-1).\]
	In the latter lower bound we can exchange $k_1$ with any other $k_i$.
\end{theorem}

In the case when $k_1=k_2=\dots =k_c=k$, Theorem~\ref{theorem:satgraphs} implies that $\sat_{\G}(k;c)\le (k-1)^c$. Using an idea of Pálvölgyi~\cite{pdperscomm}, with probabilistic methods we improve this bound in the case when $c$ is large compared to $k$.

\begin{theorem}\label{thm:random}
$\osat_{\G}(k;c)\le 48k^2c^{k^2}$.
\end{theorem}

\subsection*{Posets}

In this paper we are also interested in saturation problems on partially ordered sets (posets). Dilworth's theorem~\cite{Dilworth} answers a Ramsey-type problem about posets, since it implies that a poset of size $(k-1)(l-1)+1$ contains either a chain of length $k$ or an antichain of length~$l$. A natural saturation and semisaturation version of this problem can be posed. Let $\P$ denote the set of all finite posets. Given $P = (S, \leq_P)$, $P'=(S',\leq_{P'})\in \P$, we say $P'$ \emph{extends} $P$ if $S\subsetneq S'$ and for all $x,y\in S$, $x\leq_{P} y$ if and only if $x\leq_{P'} y$.  
A poset $P\in \P$ is $(k,l)$-semisaturated if every poset $P'\in \P$ extending $P$ contains a $k$-chain or an $l$-antichain which is not completely contained in $P$. Similarly as before, $\osat_{\P}(k,l)$ denotes the minimum size of such a semisaturated poset. If $P$ is additionally $k$-chain and $l$-antichain free, then we say that $P$ is $(k,l)$-saturated. We use $\sat_{\P}(k,l)$ to denote the minimum size of such a saturated poset. We also define $\ram_{\P}(k,l)$ as the maximum size of a $(k,l)$-saturated poset. Using this notation, Dilworth's theorem implies that 
\[\ram_{\P}(k,l)=(k-1)(l-1).\]

For the semisaturation number of posets we show the following.

\begin{theorem}\label{thm:weakGeneralPoset}
	\[\osat_{\P}(k,1)=0,\;	\osat_{\P}(k,2)=k-1.\] 
	For $l\ge 3$, we have
	\[\osat_{\P}(k,l)= \min(2k+l-5,k+3l-7).\]
\end{theorem}

For the saturation numbers of general posets, we prove the following theorem.

\begin{theorem}\label{thm:satGeneralPoset}
	\[\sat_{\P}(k,l)=(k-1)(l-1)=\ram_{\P}(k,l).\] 
\end{theorem}
When the Ramsey and the saturation numbers are the same we gain further insight into the structure of the saturated objects. For example every saturated object has the same size. Other examples of this kind include the intersecting families of subsets of an $n$ element set mentioned in the introduction and, as we will see, sequences without increasing and decreasing subsequences of given lengths.  

\subsection*{Monotone point sets and sequences}
Another well-known Ramsey-type result is the Erdős-Szekeres theorem on monotone point sets. A point set in general position is said to be monotone increasing (resp. decreasing) if when ordered according to the $x$-coordinates, the $y$-coordinates of the points are monotone increasing (resp. decreasing). The theorem of Erdős and Szekeres~\cite{ES} states that a set of $(k-1)(l-1)+1$ points in general position contains either an increasing subsequence of length $k$ or a decreasing subsequence of length $l$. There is an equivalent formulation of this result in terms of sequences. It states that a sequence of $(k-1)(l-1)+1$ numbers must contain either an increasing subsequence of length $k$ or a decreasing sequence of length $l$. We will work with both of these formulations.

We can convert this problem into a saturation problem in the usual way. A sequence  $S$ of distinct numbers (resp. point set with distinct $x$- and $y$-coordinates) is called $(k,l)$-saturated if it does not contain an increasing subsequence (resp. subset) of length $k$ or a decreasing subsequence (resp. subset) of length $l$ but any sequence (resp. point set with distinct $x$- and $y$-coordinates) $S'$ that contains $S$ as a subsequence (resp. subset) has either an increasing subsequence (resp. subset) of length $k$ or a decreasing subsequence (resp. subset) of length $l$. The functions $\sat_{\S}(k,l)$ and $\osat_{\S}(k,l)$ are defined analogously to as before. Moreover, we define $\ram_{\S}(k,l)$ to be the maximum size of a $(k,l)$-saturated sequence.
With this notation the Erdős-Szekeres theorem states that \[\ram_{\S}(k,l)=(k-1)(l-1).\]

For saturation numbers, we prove the following theorem.

\begin{theorem}\label{thm:seqSat}

\[\sat_{\S}(k,l) = (k-1)(l-1)=\ram_{\S}(k,l).\] 
\end{theorem}

In other words Theorem~\ref{thm:seqSat}  says that if a sequence of distinct numbers does not contain  an increasing subsequence of length $k$ or a decreasing subsequence of length $l$, then either we can extend the sequence without creating such a subsequence or the length of the sequence is $(k-1)(l-1)$.    

For semisaturation numbers we have the following theorem.

\begin{theorem}\label{thm:weakSatSeq}
	
	\[\osat_{\S}(1,l) = \osat_{\S}(k,1) = 0.\] 
For $k, l\geq 2$, we have 
	\[\osat_{\S}(k,l) = \min(2k+l-5, 2l+k-5).\]
\end{theorem}

\subsection*{Convex point sets}
Finally, we investigate the saturation problem for convex point sets in the plane. If a set of $n$ points is in convex position, then we say that the points form a \emph{convex $n$-set}. A set of $k$ points in convex position is called a $k$-cup (resp. $k$-cap) if the points lie on the graph of a convex (resp. concave) function (possibly multivalued).  Let $\ram_{\C\C}(k,l)$ denote the size of the largest set of points in general position which contains neither a subset forming a $k$-cup nor a subset forming an $l$-cap. Similarly, let $\sat_{\C\C}(k,l)$ denote the size of the smallest point set which contains neither a subset forming a $k$-cup nor a subset forming an $l$-cap, such that adding any new point yields a $k$-cup or $l$-cap. Finally, let $\osat_{\C\C}(k,l)$ denote the size of the smallest point set such that adding any new point introduces a new $k$-cup or a new $l$-cap.

In 1935, Erdős and Szekeres~\cite{ES} showed that \[\ram_{\C\C}(k,l)={\binom{k+l-4} {k-2}}.\]

While we were not able to obtain non-trivial bounds for the saturation problem, for the semisaturation problem we have the following result.

\begin{theorem}\label{thm:cupcapSemiSat}
	We have
	\begin{align*}
	 \osat_{\C\C}(k,3)=\osat_{\C\C}(3,k)&=k-1.\\
    \end{align*}
    For $k\geq4$, we have
	\begin{align*}
     \osat_{\C\C}(k,4)=\osat_{\C\C}(4,k)&=2k-2,
    \end{align*}
    and for $k\geq 5$ and $l \ge 5$,
	    \[2k+2l-12 \le \osat_{\C\C}(k,l)\le 2k+2l-10.\] 

\end{theorem}
The original motivation for investigating point sets free of cups and caps was to give an upper bound, namely ${ \binom{2n-4} {n-2}}$, on the maximum number of points in the plane avoiding a convex $n$-gon.  Erd\H{o}s and Szekeres provided a lower bound of size $2^{n-1}$, and after a number of subsequent improvements a nearly optimal upper bound of size $2^{n+o(n)}$ was provided by Suk~\cite{suk}.  

An intriguing problem is to obtain the analogous saturation result.

\begin{problem}\label{convex_sat}
What is the minimum possible size of a point set in the plane in general position which does not contain a convex $n$-set, but adding any extra point (in general position) creates one?    
\end{problem}

We could not even determine if the answer is polynomial in $n$ or not. Note that if we drop the general position assumption the problem becomes trivial, one can simply take $n-1$ points on a line.

For the respective semisaturation problem, let $\osat_{\C}(n)$ denote the minimum possible size of a point set in the plane general position, such that adding any extra point to it (in general position) creates a new convex $n$-set. With this notation we prove the following theorem.

\begin{theorem}\label{thm:convex}
	\[\osat_{\C}(n) = 2n-4.\]
\end{theorem}

Unlike the problem about cups and caps, this problem generalizes easily to higher dimensions. Let $\osat_{\C,d}(n)$ denote the minimum possible size of a point set in $\mathbb{R}^d$, such that it is in general position and adding one extra point to it (in general position) creates a new convex $n$-set. We obtain the following result.

\begin{theorem}\label{thm:osat00}
\[\osat_{\C,d}(n) \ge n-1 + \floor{\frac{n-2}{d}}.\]
\end{theorem}

\section{Graphs}\label{section:graph} 

\begin{proof}[Proof of Theorem \ref{theorem:satgraphs}]
	First we prove the upper bounds.
	For $c=2$, one sharp construction is when the blue (the first color) edges form the graph consisting of $l-1$ disjoint copies of $K_{k-1}$.
	Another construction is when the red (the second color) edges form the graph consisting of $k-1$ disjoint copies of $K_{l-1}$.
	It is easy to see that these two-edge-colored graphs are saturated.
	It is also easy to generalize these constructions to $c>2$ colors: Start with a graph $G_0$ with a single vertex. Now one by one for each color $i$ with $1\leq i\leq c$, construct a colored graph $G_i$ by replacing each vertex $v_j$ of $G_{i-1}$ with a clique $S_j$ of size $k_i-1$ such that all edges in the clique are in color $i$. The edges between every pair of such cliques $S_{j_1}$ and $S_{j_2}$ are in the same color as the edge $v_{j_1}v_{j_2}$ in $G_{i-1}$. It is not hard to see that the graph $G_c$ obtained by this construction is saturated.
	
	Now we prove the lower bound for $c=2$. Take a minimal two-edge-colored semisaturated graph (with respect to a blue $k$-clique and red $l$-clique). Extract maximal complete blue subgraphs (cliques) greedily one by one until we partition all the vertices into cliques (for a more advanced treatment of such greedy partitions see~\cite{gyorikeszegh}). Assume that the first $i$ blue cliques have size at least $k-1$ and the rest have size at most $k-2$. If $i\ge l-1$, then $G$ has at least $(k-1)(l-1)$ vertices and we are done. Otherwise when $i<l-1$, let $p$ be a new vertex, which we add to $G$ and connect with red edges to the first $i$ cliques and blue edges to the rest of the vertices. It is easy to see that in the resulting graph there is neither a blue clique of size $k$ containing $p$ nor a red clique of size $l$ containing $p$. Hence $G$ is not semisaturated, which contradicts our assumption.
	
	Finally, we prove the lower bound for $c>2$. Again, take a minimal $c$-edge-colored semisaturated graph $G$.
	If there is no clique of size $k_1-1$ of the first color in $G$, then we can connect a new vertex $p$ to every vertex with the first color. It follows that $G$ is not semisaturated, giving us a contradiction. Otherwise there exists a clique of size $k_1-1$ with color $1$ in $G$. Take such a clique $S$ and connect $p$ to every vertex in $S$ with the second color. Now $G-S$ must again contain a clique of size $k_1 -1$ with color $1$ otherwise we can connect $p$ to the rest of the vertices in $G$ with color one. We can connect $p$ to the vertices of this clique as well with the second color if $k_2>2$. Repeating this argument, we keep finding additional disjoint cliques of size $k_1-1$ of the first color. When we have $k_2-2$ such cliques, we continue to pull out cliques of size $k_1-1$ of the first color and connect them to $p$ with color $3$.  Continuing in this way we find altogether $(k_2-2)+(k_3-2)+\dots+(k_c-2)+1$ cliques of size $k_1-1$ of color $1$, showing that indeed the number of vertices is at least  $(k_1-1)(k_2+k_3+\dots+k_c-2c+3)$.
	As the role of the first color was not used, the same way we can find enough cliques of color $i$.
\end{proof}




In the case of two colors we have determined the exact bound for $\osat_{\G}$.  Using the following trivial observation we see that the next open case is when $c=3$ and $k_1=k$, $k_2=k_3=3$, for which Theorem~\ref{theorem:satgraphs} gives the lower bound $3(k-1)$ and upper bound $4(k-1)$ for both the saturation and semisaturation problems.
\begin{obs}\label{obs:kis2}
	\[\sat_{\G}(2,k_1,\dots, k_c)=\sat_{\G}(k_1,\dots, k_c),\]
	\[\osat_{\G}(2,k_1,\dots, k_c)=\osat_{\G}(k_1,\dots, k_c).\]
\end{obs}

We state an equivalent formulation of the $c=3$, $k_1=k$, $k_2=k_3=3$ case as an problem.  

\begin{problem}
	Is it true that the vertices of every $3$-edge colored complete graph $G$ with $n=4(k-1)$ vertices can be partitioned into three parts, the first part avoiding a $K_{k-1}$ of the first color, the second part avoiding edges of the second color and the third part avoiding edges of the third color? This is equivalent to the $\osat$ problem, in the $\sat$ variant we further assume that $G$ itself avoids $K_{k}$ of the first color and triangles of the second and third colors.
\end{problem}

We now give a probabilistic argument improving the upper bound in Theorem~\ref{theorem:satgraphs} in some cases.

\begin{proof}[Proof of Theorem \ref{thm:random}]

 Consider a uniform random coloring of the edges of the complete graph $K_{nc}$ with $c$ colors, that is, each edge is assigned one of the $c$ colors uniformly randomly and independently. We claim that the resulting edge-colored graph $G$ is semisaturated with positive probability if $n$ is large enough. 
 
 To show that $G$ is semisaturated, it suffices to show that the subgraph of $G$ induced by any set of $n$ vertices contains a monochromatic $K_{k-1}$ of each color. Indeed, if we add an extra vertex $q$ to $G$ and color the edges incident to $q$ in any way to get the graph $G'$, then by pigeonhole principle there will be $n$ edges incident to $q$ having the same color $d$. Since the endpoints of those $n$ edges (other than $q$) induce a subgraph in $G$ that contains a monochromatic $K_{k-1}$ of each color, we then can find a new $K_k$ in color $d$ in $G'$.

 Note that since we pick the color of each edge uniformly randomly and independently, each color class can be considered as an instance of the Erd\H{o}s-R\'enyi random graph $G(nc,\frac{1}{c})$. So we need to bound the probability of $G(nc,\frac{1}{c})$ having $n$ vertices whose induced subgraph does not contain a copy of $K_{k-1}$. We first need many pairwise edge-disjoint copies of $K_k$ in $K_n$. For our purposes the following simple lemma is enough:

\begin{lemma}\label{decomp}

We can find $\frac{1}{16k^2}n^2$ pairwise edge-disjoint copies of $K_k$ in $K_n$ for any $n\ge 4k^2$. 

\end{lemma}
\begin{proof}

Lemma~\ref{decomp} can be proved in many ways. For example it easily follows from the following construction. Let $\{(i,j)|i\in [k], j\in[r]\}$ be the first $kr$ vertices of the $K_n$ where $r$ is a the largest prime such that $kr\le n$. From Bertrand's postulate we know that $r$ is at least $\floor {\frac{n}{2k}}$. Since $n>4k$ we have $\floor{\frac{n}{2k}}\ge \frac{n}{2k}-1\ge \frac{n}{4k}$. Consider the cliques whose vertex set has the following form: $\{(1,a+b),(2,a+2b),\dots, (k,a+bk) \}$ where $a,b\in [r]$ and the second coordinates are understood modulo $r$. It is easy to see that this gives us $\frac{1}{16k^2}n^2$ disjoint $K_k$'s. Indeed, suppose otherwise that two copies share an edge, then for some values $(a,b)\ne(a',b')$, $i\ne j$ we would have $a+ib=a'+ib'$ and $a+jb=a'+jb'$ modulo $r$, implying $(i-j)(b-b')=0$ modulo $r$. As $1\le i,j$ and $i,j\le k\le r$ ($n\ge 4k^2$ implies that $k\le r$), this gives $b=b'$, which in turn implies $a=a'$, a contradiction.
\end{proof}

Now we can calculate a bound on the probability of $G(n,p)$ containing a $K_{k-1}$ (where $p=1/c$). Let $n\ge 4k^2$ to be chosen later. First we fix $\frac{1}{16k^2}n^2$ disjoint copies of $K_{k-1}$ using Lemma \ref{decomp}. For each copy the probability that it is not in $G(n,p)$ is $1-p^{\binom{k-1}{2}}$. Since the cliques are edge-disjoint, the probability that no $K_{k-1}$ appears at all is at most $\lp 1-p^{\binom{k-1}{2}}\rp^{\frac{n^2}{16k^2}}\le e^{-p^{\binom{k-1}{2}}\frac{n^2}{16k^2}}$.

Returning to the original problem we see that there are $\binom{cn}{n}\le (ec)^n$ ways to choose $n$ vertices out of $nc $. The colors have a symmetric role in the problem. Hence by the union bound, the probability that we can find $n$ vertices and a color such that there is no $K_{k-1}$ of that color among those $n$ vertices is at most \[c(ec)^n e^{-p^{\binom{k-1}{2}}\frac{n^2}{16k^2}}.\]

Picking $n=3\log(c)16k^2c^{\binom{k-1}{2}}$ we get 

\[c(ec)^n e^{-p^{\binom{k-1}{2}}\frac{n^2}{16k^2}}\le e^nc^{n+1}e^{-3\log(c)n}< 1.\]

Therefore the probability of the bad cases is less than $1$. So we can find a semisaturated graph on $cn=c\cdot 3\log(c)16k^2c^{\binom{k-1}{2}}\le 48k^2 c^{k^2}$ vertices. 
\end{proof}

\section{Posets}\label{section:poset}

\begin{proof}[Proof of Theorem \ref{thm:weakGeneralPoset}]
	If $l=1$, then obviously adding an element to the empty poset will introduce an antichain of size $1$. Thus $\osat_{\P}(k,1)=0$.
	
Consider the case when $l=2$. If a newly added element is incomparable with any of the elements of the poset, then we find a new antichain of size two. So if $P$ is a poset which is not semisaturated for $l=2$, then we must be able to add an element comparable to all elements of $P$ without introducing a new chain of size $k$. This is possible if and only if $P$ does not contain a chain of size $k-1$. On the other hand, the smallest poset containing a chain of size $k-1$ is the poset containing only this chain on $k-1$ elements, and this poset is clearly semisaturated for $2$-antichains and $k$-chains. Thus $\osat_{\P}(k,2)=k-1$ (and this is the only semisaturating poset of this size).

Now we may assume that $l\ge 3$. We show two semisaturated posets, one with $2k+l-5$ elements and one with $k+3l-7$ elements. 

For the first construction consider an antichain $A$ of size $l-1$ and then add two chains, $C_1$ and $C_2$, of length $k-2$ to the poset such that the $C_1$ lies below the elements of $A$ and $C_2$ lies above them (see Figure~\ref{fig:theo5}). The resulting poset is semisaturated. Indeed, if we add $p$ such that it is not comparable to any element of $A$, then $A\cup\{p\}$ is an antichain of length $l$. If $p$ lies below an element $a$ in $A$, then $C_2\cup\{a,p\}$ is a chain of length $k$. Similarly if $p$ lies above an element $a$ in $A$ then $C_1\cup\{a,p\}$ is a chain of length $k$. Therefore we cannot add $p$ without creating a chain of size $k$ or an antichain of size $l$. 

\begin{figure}[!ht]
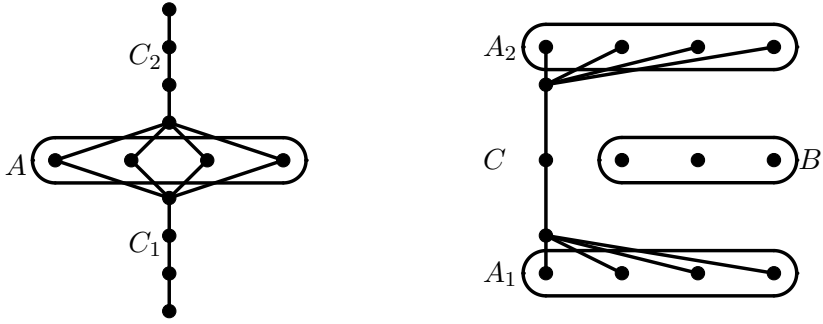

    \centering
    \begin{minipage}{0.4\textwidth}
    \include{posetconst4}
    \end{minipage}
    \begin{minipage}{0.4\textwidth}
    \include{posetconst3}
    \end{minipage}
    \caption{The two constructions for $k=6$, $l=5$.}
    \label{fig:theo5}
\end{figure}

The second construction starts with two disjoint antichains $A_1$ and $A_2$ of length $l-1$. Then we add a chain $C$ of length $k-3$ between them, that is, every element of $C$ is above every element of $A_1$ and below every element of $A_2$. Finally we add an antichain $B$ of $l-2$ elements such that they are incomparable to everything.  

To see that this construction is semisaturated suppose that we can add an element $p$ without creating a chain of size $k$. Then $p$ cannot be above any element of $A_2$ nor below any element of $A_1$. On the other hand $p$ must be comparable to some elements $a_1\in A_1$ and $a_2\in A_2$ otherwise we would get an antichain of length $l$. Consequently, $p$ is above $a_1$ and below $a_2$. We know that $C\cup \{a_1,a_2\}$ is a chain of length $k-1$ so there must be an element $c\in C$ such that $p$ and $c$ are incomparable otherwise $C\cup\{a_1,a_2,p\}$ would be a chain of size $k$. Since $B\cup \{c\}$ is an antichain and $B \cup \{p,c\}$ cannot be an antichain, $p$ must be comparable to some element $q\in B$. If $p$ is above $q$ then $a_2$ is comparable to $q$ through $p$. If $p$ is below $q$ then $a_1$ is comparable to $q$ through $p$. But neither case is possible since $q$ is incomparable to both $a_1$ and $a_2$ in the original poset.

To show that we need at least $\min(2k+l-5,k+3l-7)$ elements for semisaturation we start with the following observation.  Let $P$ be a semisaturated poset, and let $L$ denote those elements that are the minimal elements of a chain of length $k-1$ in $P$. We claim that $|L|\ge l-1$. To see this, observe that if we add an element $p$ below every element of $P\setminus L$ and incomparable to every element of $L$ (this is possible since $L$ is clearly a downset in the poset), then no chain of length $k$ is created. Since $P$ is semisaturated, it follows that $p$ must be in an antichain of length $l$ which lies in $L\cup \{p\}$. Thus, $|L|\ge l-1$ and $L$ contains an antichain $L'$ of size $l-1$. Similarly we define $U$ to be the set of elements that are maximal elements of a chain of length $k-1$. In the same way we can see that $|U|\ge l-1$, and $U$ contains an antichain $U'$ of size $l-1$.

If $U\cap L\ne \emptyset$, then there is a chain of length $2k-3$. Since this chain intersects $L'$ in at most one element, we have at least $2k-3+l-1-1=2k+l-5$ elements. 

If $U\cap L= \emptyset$ and $1\le k\le 3$, then the number of elements is at least $|U|+|L|\ge 2l-2\ge 2k+l-5$, as required (using that $l\ge 3$). From now on we assume $k\ge 3$ and consider two cases. First suppose that every element of $U$ is comparable to every element of $L$. Then we can add $p$ to the poset such that it lies below every element of $U$, above every element of $L$ and incomparable to every other element. 

Suppose $p$ creates a new chain $C$ of length $k$. Clearly $C\subset U\cup L\cup \{p\}$. By symmetry we may assume that $|C\cap U|\ge \ceil{\frac{k-1}{2}}$. Let $u$ be the first element in $C$ above $p$. Since $u\in U$ there is a chain $C_2$ of size $k-1$ ending in $u$. Consequently $(C\cap U) \cup C_2$ is a chain of length at least $k-1+\ceil{\frac{k-1}{2}}-1$. In total we have found $|U'\cup L'\cup ((C\cap U) \cup C_2)| $ elements. Any chain intersects any antichain in at most one element, therefore \[|U'\cup L'\cup ((C\cap U) \cup C_2)|\ge l-1+l-1+k-1+\ceil{\frac{k-1}{2}}-1-2=2l+k+\ceil{\frac{k-1}{2}}-6.\] As $\min(k+3l-7,2k+l-5)$ is at most \[\floor{\frac{(2k+l-5)+(k+3l-7)}{2}}=\floor{2l+\frac{3}{2}k-6}=2l+k+\floor{\frac{k}{2}}-6=2l+k+\ceil{\frac{k-1}{2}}-6,\] we have at least $\min(k+3l-7,2k+l-5)$ elements.

Suppose now that adding $p$ does not create a new chain of length $k$. Then since $P$ is $(k,l)$-semisaturated, it must happen that adding $p$ creates a new antichain of length $l$. Since $p$ is comparable to the elements of $U$ and $L$ we have found $l-1$ new elements that are not in $L\cup U$. Therefore we have three disjoint antichains of length $l-1$ ($L'$, $U'$ and this antichain containing $p$). We must also have a chain of length $k-1$ and this chain intersects each of the three antichains in at most one element. Therefore we have at least $3(l-1)+k-1-3=k+3l-7$ elements.   

The only remaining case is when  $U\cap L= \emptyset$ and we can find $u\in U$ and $q\in L$ such that $u$ and $q$ are incomparable. This implies that the chains going up from $q$ and going down from $u$ are disjoint. So we have two disjoint chains of length $k-1$, giving us at least $2(k-1)+2(l-1)-4=2k+2l-8\ge 2k+l-5$ elements.
\end{proof}	
	
\begin{proof}[Proof of Theorem \ref{thm:satGeneralPoset}]
	Given a poset $P$ with fewer than $(k-1)(l-1)$ elements which contains no $k$-chain and no $l$-antichain, we need to show that we can always add an element $p$ to $P$ in such a way that the resulting poset still avoids $k$-chains and $l$-antichains. 

	If the maximum length of an antichain in $P$ is at most $l-2$, then we can easily add a new element to $P$ incomparable with all elements of $P$, and thus still avoid $k$-chains and $l$-antichains.

	Suppose now that the size of the maximal antichain is $l-1$. By Dilworth's theorem we can decompose $P$ into $l-1$ chains. By our assumption that we are $k$-chain free, all of them have size at most $k-1$ and at least one of them has size strictly less than $k-1$. Denote one such chain by $C$. For an element $c\in C$ denote by $D_c$ the subchain of $C$ consisting of $c$ and the elements that are below $c$ in $C$. Similarly, $U_c$ is the subchain of $C$ containing $c$ and the elements above $c$ in $C$.
	
	First suppose that there is no chain of size $k-1$ in $P$ whose bottom element is the bottom element $q'$ of $C$. Then we add a new element $p$ directly under $q'$ and incomparable with all the elements that are not above $q'$ (Figure~\ref{fig:theo6a}, left side). We claim that $P\cup \{p\}$ still avoids $k$-chains and $l$-antichains and so $P$ was not saturated, a contradiction. Indeed, as the poset can still be partitioned into $l-1$ chains (the former partition of $P$ with the difference that $C$ is extended with $p$ as a new bottom element under $q'$), it follows that there is no antichain of length $l$. Also a chain of length $k$ must have $p$ as its bottom element, and then $q'$ as its element directly above $p$, but then this chain minus $p$ would be a chain of length $k-1$ in $P$ with bottom element $q'$, contradicting our assumption. The case when there is no chain of size $k-1$ in $P$ whose top element is the top element of $C$ is handled similarly.
\begin{figure}[!ht]
    \centering\begin{minipage}{0.4\textwidth}
    \includegraphics{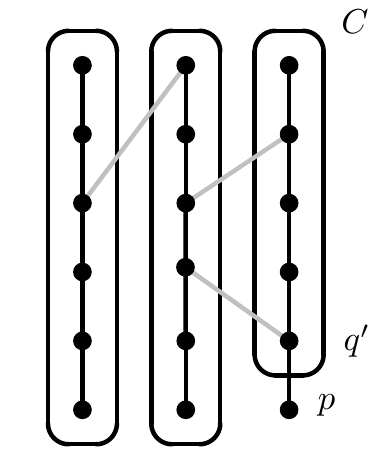}
\end{minipage}
    \begin{minipage}{0.4\textwidth}
    \includegraphics{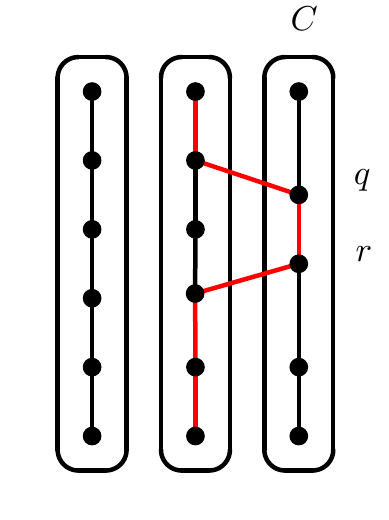}
\end{minipage}
    \caption{Adding $p$ to the bottom of $C$ and finding a long chain if $q$ is above $r$.}
    \label{fig:theo6a}
\end{figure}

	Thus, we may assume that there is a largest element $q$ of $C$ for which there is a chain $C_q$ of size $k-1$ containing $q$ whose part below $q$ (including $q$) coincides with $D_q$. Similarly, there is a smallest element $r$ of $C$ for which there is a chain $C_r$ of size $k-1$ containing $r$ whose part above $r$ (including $r$) coincides with $U_r$. 

We claim that $r$ is above $q$. Suppose on the contrary that $q$ is above $r$ or that they coincide. Then taking the part of $C_q$ above $q$, the part of $C$ between $q$ and $r$ and the part of $C_r$ below $r$ we get a chain $C'$ whose length is at least $k$, a contradiction. Indeed, the sum of $\abs{C}$ and $\abs{C'}$ is the same as the sum of $\abs{C_q}$ and $\abs{C_r}$, and so as $C_q$ and $C_r$ have $k-1$ elements and $C$ has at most $k-2$ elements, it follows that $C'$ must have at least $k$ elements, a contradiction (see Figure~\ref{fig:theo6a}).
	
Thus $q$ is not the top element of $C$ and there is an element $s$ of $C$ directly above it. Now add a new element $p$ directly above $q$ and below $s$ such that $p$ is incomparable to all elements that are not below $q$ or above $s$ (Figure~\ref{fig:theo6b}). We claim that $P\cup \{p\}$ still avoids $k$-chains and $l$-antichains and thus $P$ was not saturated, a contradiction. The poset $P\cup \{p\}$ can still be partitioned into $l-1$ chains by taking the previous partition except that $C$ is extended with $p$ put between $q$ and $s$. Thus, the resulting poset still avoids $l$-antichains. Suppose now that it contains a $k$-chain, it necessarily contains $p$, and then $q$ is directly below $p$ in the chain and $s$ is directly above $p$ in the chain. Deleting $p$ from this $k$-chain we get a $(k-1)$-chain $C'$ in $P$ which has $q$ and $s$ directly above each other. 
	
Let $D'_q$ be the part of $C'$ below $q$ (including $q$) and $U'_{s}$ be the rest of $C'$, that is the part of $C'$ above $s$ (including $s$). By the definition of $q$ there is a chain of size $k-1$ whose bottom part is $D_q$, thus $D'_q$ can have size at most as much as $D_q$ otherwise the top part of this chain extended with $D'_q$ would be a chain of size at least $k$. Again by the maximality of $q$, the chain formed by $D_q$ and $U'_{s}$ can have size at most $k-2$, thus the chain formed by $D'_q$ and $U'_{s}$ can also have size at most $k-2$, but this is exactly $C'$ which was of size $k-1$, a contradiction.
\begin{figure}[!ht]
    \centering
    \centering\begin{minipage}{0.4\textwidth}
        \includegraphics[width=5cm]{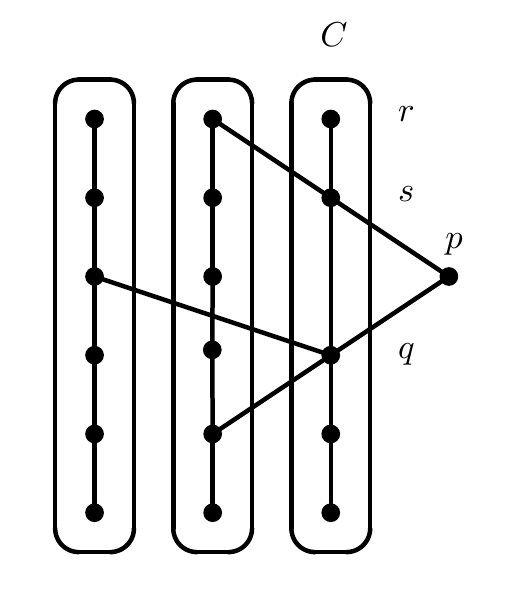}
\end{minipage}
    \begin{minipage}{0.4\textwidth}
    \includegraphics[width=5cm]{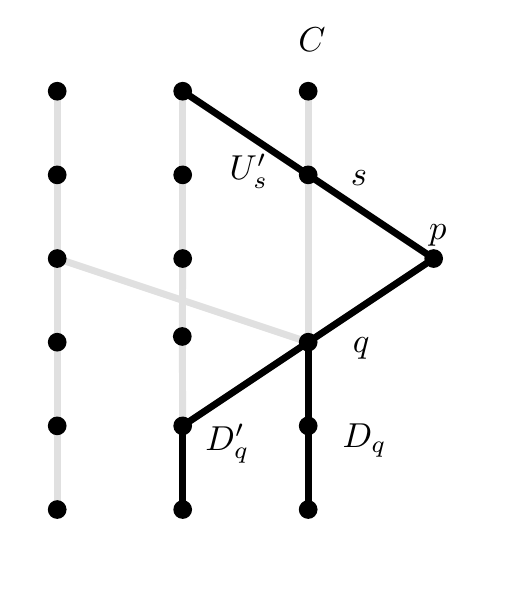}
\end{minipage}
    \caption{Inserting $p$ into the poset.}
    \label{fig:theo6b}
\end{figure} 
\end{proof}

\section{Saturation of monotone point sets and sequences}\label{section:strongseq}
Throughout this section, the term sequence will always refer to a sequence of distinct real numbers. We may, without loss of generality, assume that the numbers in the sequence are positive.

\begin{definition}
A sequence $S$ is $(k,l)$-saturated if $S$ contains no increasing subsequence of length $k$ nor decreasing subsequence of length $l$, but any sequence $S'$ containing $S$ as a proper subsequence contains either an increasing subsequence of length $k$ or a decreasing subsequence of length $l$. Let $\sat_{\S}(k,l)$ denote the minimum possible length of a $(k,l)$-saturated sequence. 
\end{definition}

Let $\vec{a} = [a_1, a_2,\dots, a_m]$ be a $(k+1,l+1)$-saturated sequence. We define the function ${\gamma_{\vec{a}}:[m] \to [l] \times [k]}$ by $\gamma_{\vec{a}}(t)=(i,j)$, where $i$ is the length of the longest decreasing subsequence of $\vec{a}$ ending at $a_t$ and $j$ is the length of the longest increasing subsequence of $\vec{a}$ ending at $a_t$. Since for every $n' < n\le m$ we can extend either the longest increasing subsequence or the longest decreasing subsequence ending at $a_{n'}$ by appending $a_n$ to the end, we have the following observation.   

\begin{obs}\label{obs:inc}
    If $n'<n$, then at least one coordinate of  $\gamma_{\vec{a}}(n)$ is strictly larger than the corresponding coordinate of $\gamma_{\vec{a}}(n')$. 
\end{obs}
 Define an $l \times k$ matrix $R_{\vec{a}} = (r_{ij})$ corresponding to the sequence $\vec{a}$ by setting
\[r_{ij} = \begin{cases}
			\gamma_{\vec{a}}^{-1}(i,j) & \textrm{if $(i,j) \in \image(\gamma_{\vec{a}})$,} \\
            0 & \textrm{otherwise.}
	   \end{cases}\]

Define an $l \times k$ matrix $V_{\vec{a}} = (v_{ij})$ corresponding to the sequence $\vec{a}$ by setting 
\[v_{ij}=
	  \begin{cases}
	a_{\gamma_{\vec{a}}^{-1}(i,j)} & \textrm{if $(i,j) \in \image(\gamma_{\vec{a}})$,} \\
    0 & \textrm{otherwise.}
      \end{cases}
    \]
    
Finally, let $W_{\vec{a}}=(w_{ij})$ be the $l \times k$ matrix such that the $i$-th row of $W_{\vec{a}}$ is the $(l+1 -i)$-th row of $V_{\vec{a}}$ for every $i \in [l]$. For example, consider the sequence $\vec{a}=[33,11,22,55,44]$ with $k=3$, $l=2$. In this case $R_{\vec{a}}=\begin{pmatrix}
1 & 0 & 4\\
2 & 3 & 5
\end{pmatrix}$ and $V_{\vec{a}}=\begin{pmatrix}
33 & 0 & 55\\
11 & 22 & 44
\end{pmatrix}$ and thus $W_{\vec{a}}=\begin{pmatrix}11 & 22 & 44\\
33 & 0 & 55
\end{pmatrix}$.

\begin{definition}
We say that a matrix $(m_{ij})$ is \emph{partially increasing} if all the elements are nonnegative, the positive values are distinct and $i_1\le i_2$, $j_1\le j_2$ implies $m_{i_1j_1}\le m_{i_2j_2}$ for all nonzero elements of the matrix.
\end{definition}

\begin{lemma}\label{lem:young_tableau}
$R_{\vec{a}}$ and $W_{\vec{a}}$ are both partially increasing. 
\end{lemma}

\begin{proof}
	Observation~\ref{obs:inc} implies that $R_{\vec{a}}$ is partially increasing.
	
	To show that $W_{\vec{a}}$ is partially increasing we need to prove that if we take $i_1\ge i_2$ and $j_1\le j_2$, then $v_{i_1,j_1}\le v_{i_2,j_2}$ whenever $v_{i_1,j_1}$ and $v_{i_2,j_2}$ are positive numbers. Assume on the contrary that $v_{i_1,j_1}> v_{i_2,j_2}>0$. Let $n, n'$ be the indices such that $v_{i_1,j_1}=a_n>a_{n'}=v_{i_2,j_2}$. 	
	
	First, if $n<n'$, then a decreasing subsequence of length $i_1$ ending in $a_{n}$ can be extended with $a_{n'}$ and thus the longest increasing subsequence ending in $a_{n'}$ has length at least $i_1+1$, which contradicts that $i_2 \leq i_1$.
	
	Second, if $n>n'$, then an increasing subsequence of length $j_2$ ending in $a_{n'}$ can be extended with $a_{n}$ and thus the longest increasing subsequence ending in $a_{n}$ has length at least $j_2+1$, which contradicts that $j_1 \leq j_2$.
\end{proof}

Let $\mathbb{S}_{k,l}$ be the set of extremal $(k+1, l+1)$-saturated sequences of length $kl$ whose entries are distinct integers in $[kl]$. Observe that when $\vec{a} \in \mathbb{S}_{k,l}$, $R_{\vec{a}}$ and $W_{\vec{a}}$ have all positive entries, and are increasing in both rows and columns by Lemma~\ref{lem:young_tableau}.

Before moving on to the proof of Theorem~\ref{thm:seqSat}, we briefly discuss how our results relate to the classification of extremal sequences for the Erd\H{o}s-Szekeres theorem in terms of Young tableaus. As $R_{\vec{a}}$ and $W_{\vec{a}}$ have values from $[kl]$, they correspond to a pair of standard rectangular Young tableaus $\mathbb{Y}_{l,k} \times \mathbb{Y}_{l,k}$ (with entries in $[kl]$).
It was observed earlier by Knuth [\cite{Knuth}, Exercise 5.1.4.9] (see also [\cite{Stanley}, Example 7.23.19(b)]) that the set $\mathbb{S}_{k,l}$ is in bijection with the set of pairs of standard Young tableaus $\mathbb{Y}_{l,k} \times \mathbb{Y}_{l,k}$ (with entries in $[kl]$) via the Robinson-Schensted correspondence. Romik~\cite{romik} (see also~\cite{Czabarka-W}) gave an explicit bijection via the function $\phi(\vec{a}) = (R_{\vec{a}},W_{\vec{a}})$. 

Theorem~\ref{thm:seqSat} shows that all $(k+1,l+1)$-saturated sequences are in fact extremal, i.e., have length $kl$. Hence there is also a bijection between the set of all $(k+1,l+1)$-saturated sequences and the set of pairs of standard rectangular Young tableaus (with entries in $[kl]$).

\begin{proof}[Proof of Theorem~\ref{thm:seqSat}]
	
Let $\vec{a} = [a_1, a_2,\dots, a_m]$ be a $(k+1,l+1)$-saturated sequence. By definition $\sat_{\S}(k+1,l+1) \leq m \leq \ram_{\S}(k+1,l+1)=kl$. To simplify our notation let $\gamma=\gamma_{\vec{a}}$, $R=R_{\vec{a}}$, $V=V_{\vec{a}}$ and $W=W_{\vec{a}}$.

By the Erdős-Szekeres theorem for sequences, $\sat_{\S}(k+1,l+1) \leq k l$. Hence it suffices to show that $\sat_{\S}(k+1,l+1) \geq kl$. Let $\vec{a} = [a_1, a_2,\dots, a_m]$ be an arbitrary sequence of length $m < kl$ containing no increasing subsequence of length $k+1$ and no decreasing subsequence of length $l+1$. We need to show that $\vec{a}$ is not saturated.


\begin{claim}\label{cl:completion}  If a partially increasing matrix $M$ contains a $0$, then that $0$ can be replaced by a positive number so that the resulting matrix is still partially increasing.  
\end{claim}

\begin{proof}[Proof of Claim~\ref{cl:completion}]

Let $(i_0,j_0)$ be the position of a $0$ in $M$. If there is no nonzero entry in any position $(i,j)$ with $i_0 \le i$, $j_0 \le j$, then replace 0 with any number larger than all entries of the matrix. Otherwise let $t=\min\limits_{i\ge i_0,j\ge j_0}(M)_{ij}$. Change the value of $M$ at $(i_0,j_0)$ to be $t-\epsilon$ where $\epsilon$ is smaller than the difference of any two nonzero values of $M$ and it is also smaller than $t$. It is easy to see that since $M$ is partially increasing, the new matrix obtained is also partially increasing.     
\end{proof}

By Claim~\ref{cl:completion}, if $r_{i,j}=0$ (and thus $v_{i,j}=w_{l+1-i,j}=0$), we can replace these $0$'s with some positive number (not necessarily integers) such that $R$ and $W$ are still partially increasing. We then relabel the elements of $R$ (if necessary) with an initial sequence of integers in $[k l]$ while respecting their order. Call the resulting matrices $R'$ and $W'$ respectively, we get $V'$ from $W'$ by reversing the order of its rows, that is the $i$-th row of $V'$ is the $(l+1-i)$-th row of $W'$. Continuing the example above we obtain  that ${R'=\begin{pmatrix}
1 & 3 & 5\\
2 & 4 & 6
\end{pmatrix}}$, $V'=\begin{pmatrix}
33 & 50 & 55\\
11 & 22 & 44
\end{pmatrix}$ and $W'=\begin{pmatrix}
11 & 22 & 44\\
33 & 50 & 55
\end{pmatrix}$.

Now we can extend $\vec{a}$ as follows: Insert the number $(V')_{i,j}$ immediately after the {$((R')_{i,j}-1)${-th}} position of $\vec{a}$. Call the resulting sequence $\vec{b}$ (in our example $\vec{b}=[33,11,50,22,55,44]$). We see that $V'$ records the values in $\vec{b}$ and $R'$ records the order of these values. We want to show that $\vec{b}$ contains no increasing subsequence of length $k+1$ nor a decreasing subsequence of length $l+1$. 

Suppose $\vec{b}$ contains an increasing subsequence of length $k+1$. Then at least two elements of this subsequence must be in the same column of $V'$ by the pigeonhole principle.
Since $V'$ is just $W'$ reversed, it is decreasing in the columns. Hence, the lower one of these two elements is smaller, and since they are in an increasing subsequence it must appear earlier in $\vec{b}$. On the other hand $R'$ is increasing in the columns so the lower must come later, a contradiction. Similarly, $\vec{b}$ does not have a decreasing subsequence of length $l+1$. Hence $\vec{a}$ is not saturated and the proof is complete.
\end{proof}

\section{Semisaturation of monotone point sets and sequences}\label{section:seq}

In this section we use the point set formulation of the problem, the term point set always refers to a set of points in general position (i.e., no two points share a common $x$ or $y$ coordinate). We start with the following trivial lemma.
\begin{lemma}
	\label{lem:intersect}
	If $I$ is an increasing subset and $D$ is a decreasing subset of the point set $P$, then they intersect in at most one element.
\end{lemma}

\begin{proof}[Proof of Theorem \ref{thm:weakSatSeq}]
	Fix some $n\in\mathbb{Z}^{+}$, and assume that $\pointSet$ is semisaturated with respect to monotone
	$n$-sequences.  Then we know that any point not already contained in
	$\pointSet$ must combine with $n-1$ points from $\pointSet$ to form a
	monotone $n$-sequence.  Note that any subsequence of a monotone
	sequence must itself be a monotone sequence; as such, our analysis
	will focus on monotone $(n-1)$-sequences in $\pointSet$, and consider
	which points in the plane can be added to a given $(n-1)$-sequence to
	produce a monotone $n$-sequence.  We say an $(n-1)$-sequence {\em
		blocks} such points; thus, we can say that $\pointSet$ is saturated
	with respect to monotone $n$-sequences if and only if every point in
	the plane is either contained in $\pointSet$, or blocked by some
	monotone $(n-1)$-sequence from $\pointSet$.
	
	Consider some increasing $(k-1)$-sequence
	$\point[1],\dots,\point[k-1]$.  The set of points blocked by the
	sequence is precisely the union $\cup_{i=1}^{k}\region[i]$ of regions given by
	
	\begin{minipage}{0.55\textwidth}
		\begin{alignat*}{5}
		&\region[1] &&=& (-\infty,x_1]&\times(-\infty,y_1];\\
		&\region[i+1] &&=& [x_i,x_{i+1}]&\times[y_i,y_{i+1}],\text{ for $i=1,\dots,k-2$;}\\
		&\region[k] &&=& [x_{k-1},\infty)&\times[y_{k-1},\infty).
		\end{alignat*}
	\end{minipage}
	\begin{minipage}{0.4\textwidth}
		\includegraphics{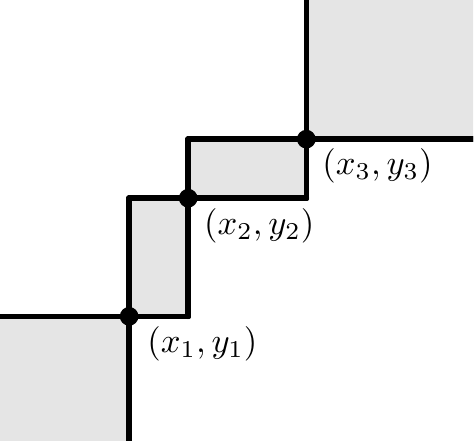}
	\end{minipage}
	\smallskip
	
	Decreasing sequences behave similarly. For our proof, we focus on points that are outside some fixed axis-parallel rectangle that contains all the points in $\pointSet$, and so are only interested in regions
	$\region[1]$ and $\region[k]$ above. More precisely, we pick bounding values $\lo{x}$, $\hi{x}$, $\lo{y}$, and $\hi{y}$ such that for each point $\coords\in\pointSet$ we have that	$\lo{x}<x<\hi{x}$ and $\lo{y}<y<\hi{y}$.  We will focus on how points
	lying on the lines forming the boundary of this region are blocked;
	note that these inequalities guarantee that such points can {\em only}
	be blocked by being at one end or the other of an increasing $(k-1)$-sequence or decreasing $(l-1)$-sequence.

	To begin, consider points along the line $y=\lo{y}$; fix any such
	point $(x,\lo{y})$.  Now, since $\lo{y}$ is strictly less than the
	$y$-coordinate of any point in $\pointSet$, we may conclude that if
	any $(n-1)$-sequence of points $\point[1],\dots,\point[n-1]$ from
	$\pointSet$ block $(x,\lo{y})$, it must be the case either that the
	sequence is decreasing (and $n=l$) and $x_{n-1}\le x$ or that the sequence is
	increasing (and $n=k$) and $x\le x_{1}$.  Viewed from the opposite perspective, we
	can see that any decreasing $(l-1)$-sequence
	$\point[1],\dots,\point[l-1]$ in $\pointSet$ blocks the left-bounded
	interval $[x_{l-1},\infty)$ on the line $y=\lo{y}$.  Symmetrically,
	any increasing $(k-1)$-sequence blocks the right-bounded interval
	$(-\infty,x_1]$.  For the entire line $y=\lo{y}$ to be blocked, then,
	we need a left-bounded interval and a right-bounded interval that
	intersect each other; this equates to a decreasing $(l-1)$-sequence
	and an increasing $(k-1)$-sequence such that the former lies entirely
	to the left of the latter.
	
	Let $\lo{D}$ and $\lo{I}$ be a decreasing $(l-1)$-sequence and increasing
	$(k-1)$-sequence from $\pointSet$, respectively, such that for all
	$\coords\in\lo{D}$ and all $\coordsP\in\lo{I}$, we have that $x \le
	x'$.  The preceding argument guarantees the existence of such, and
	furthermore a symmetric argument with respect to the line $y=\hi{y}$
	gives us that we can find an increasing $(k-1)$-sequence and a decreasing $(l-1)$-sequence
	$\hi{I}$ and $\hi{D}$, respectively, in $\pointSet$, such that for all
	$\coords\in\hi{I}$ and all $\coordsP\in\hi{D}$ we have that $x \le x'$.
	
	Our claim is that $\abs{\lo{D}\cup\lo{I}\cup\hi{D}\cup\hi{I}}\ge \min(2k+l-5,2l+k-5)$.
	We break our analysis into three cases.
	\begin{itemize}
		\item[Case:] $\lo{D}\cap\hi{D}=\emptyset$.  Recall that we have assumed
		that no two points in $\pointSet$ share either a common $x$-value or
		a common $y$-value; thus, Lemma~\ref{lem:intersect} tells us that
		$\abs{\lo{I}\cap\lo{D}}\le1$ and $\abs{\lo{I}\cap\hi{D}}\le1$.
		Thus, we get that
		\begin{align*}
		\abs{\lo{D}\cup\lo{I}\cup\hi{D}\cup\hi{I}}
		&\ge \abs{\lo{I}\cup\lo{D}\cup\hi{D}}\\
		&\ge \abs{\lo{I}}+\abs{\lo{D}}+\abs{\hi{D}}-\abs{\lo{I}\cap\lo{D}}-\abs{\lo{I}\cap\hi{D}}-\abs{\lo{D}\cap\hi{D}}\\
		&\ge k-1 + 2(l-1) - 2\\
		&= 2l+k -5,
		\end{align*}
		exactly as desired.
		\item[Case:] $\lo{I}\cap\hi{I}=\emptyset$.  This case is symmetric to
		the preceding one.  Applying Lemma~\ref{lem:intersect} appropriately gives us
		that
		\begin{equation*}
		\abs{\lo{D}\cup\lo{I}\cup\hi{D}\cup\hi{I}}
		\ge \abs{\lo{D}\cup\lo{I}\cup\hi{I}}
		\ge 2k + l -5,
		\end{equation*}
		once again.
		\item[Case:] $\abs{\lo{D}\cap\hi{D}},\abs{\lo{I}\cap\hi{I}}>0$.  Let
		$\coords\in\lo{D}\cap\hi{D}$.  Now, by our definitions of $\hi{I}$
		and $\hi{D}$, we must have that for all $\coordsP\in\hi{I}$, $x' \le
		x$ holds.  Similarly, our definitions of $\lo{I}$ and $\lo{D}$
		ensure that for all $\coordsP\in\lo{I}$, we have $x \le x'$.
		Consider combining this with our assumption that no two points in
		$\pointSet$ share either a common $x$-value or a common $y$-value.
		Say $\hi{I}$ consists of the point sequence
		$\point[1],\dots,\point[k-1]$, and $\lo{I}$ consists of
		$\point[1]',\dots,\point[k-1]'$.  Then we must have that
		\begin{equation*}
		x_1 < x_2 < \dots < x_{k-1} \le x \le x_1' < x_2' <\dots < x_{k-1}'.
		\end{equation*}
		By assumption, however, we have that
		$\lo{I}\cap\hi{I}\neq\emptyset$; this can only hold, then, if
		$p_{k-1}=p_{1}'$.  Thus, we have that $\lo{I}\cup\hi{I}$ is, in
		fact, an increasing $(2k-3)$-sequence.  A symmetric argument implies
		that, similarly, $\lo{D}\cup\hi{D}$ is a decreasing
		$(2l-3)$-sequence.  So applying Lemma~\ref{lem:intersect} gives us that
		\begin{align*}
		\abs{\lo{D}\cup\lo{I}\cup\hi{D}\cup\hi{I}}
		& \ge 
		\abs{\lo{D}\cup\hi{D}} + \abs{\lo{I}\cup\hi{I}}
		-\abs{(\lo{D}\cup\hi{D})\cap(\lo{I}\cup\hi{I})}
		\\ &    \ge (2l-3)+(2k-3) - 1
		\ge \min(2k+l-5,2l+k-5),
		\end{align*}
	\end{itemize}
	since $k$ and $l$ are at least $2$.
	In every case, the above gives us the desired lower bound, namely
	that $\abs{\pointSet}\ge \min(2k+l-5,2l+k-5)$.
	
	The upper bound follows from the following simple construction (see Figure~\ref{fig:seq}). 
	
	\begin{figure}[!ht]
		\centering
		\includegraphics{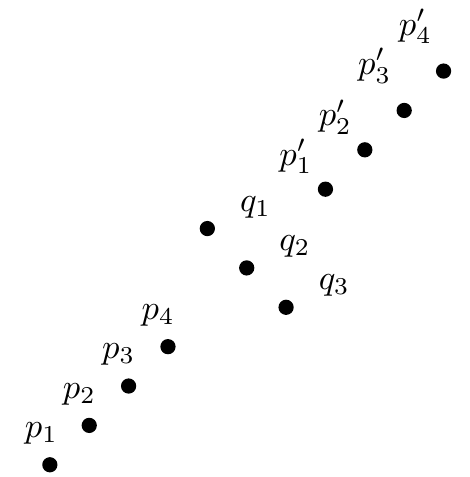}
		\caption{Construction for semisaturated sequences ($k=6,l=4$).}
		\label{fig:seq}
	\end{figure}
	
	\begin{construction}
		We present a construction of size $2k+l-5$, a construction of size $2l+k-5$ is attained by taking this construction for $k'=l$ and $l'=k$ of size $2k'+l'-5=k+2l-5$ and then reversing the order of the $x$-coordinates. Take an increasing $(k-2)$-sequence $p_1,p_2,\dots,p_{k-2}$ and let $p_{k-2}=(a,b)$.  Take another increasing $(k-2)$-sequence $p_1',p_2',\dots,p_{k-2}'$, let $p'_1 = (a',b')$ and assume $a<a'$ and $b<b'$. 
		Consider the rectangle defined by the vertices  $(a,b),(a,b'),(a',b),(a',b')$. Let $0<\epsilon< \min\left((a'-a)/4,(b'-b)/4\right)$ and consider the rectangle with corners $(a+\epsilon,b+\epsilon),(a+\epsilon,b'-\epsilon),(a'-\epsilon,b+\epsilon),(a'-\epsilon,b'-\epsilon)$. Take a decreasing $(l-1)$-sequence $q_1,q_2,\dots,q_{l-1}$ with $q_1 = (a+\epsilon,b'-\epsilon)$ and $q_{l-1} = (a'-\epsilon,b+\epsilon)$. \qedhere
	\end{construction}
\end{proof}

\section{Convex point sets}\label{section:conv}

Given a point set $S \subseteq \mathbb{R}^d$, we use $\conv(S)$ to denote the {\em convex hull} of $S$, which is the smallest convex set in $\mathbb{R}^d$ that contains $S$, and we use $\interior(S)$ to denote the interior of $S$.

First we prove the semisaturation result about cups and caps.

\begin{proof}[Proof of Theorem \ref{thm:cupcapSemiSat}]	
For $k\le 2$ and $l\le 2$ the problem is trivial.  Any two points form a 2-cup and also a 2-cap. Hence $\osat_{\C\C}(2,l)=\osat_{\C\C}(k,2)=1$. From now on let us assume that $k\ge 3$ and $l\ge 3$. 

Let $P$ be a point set that is semisaturated for $3$-cups and $l$-caps. Let $L$ be the set of lines determined by the points of $P$. There is an unbounded region $R$ in the plane bounded by parts of the lines that lies below every line of $L$. If we add a point $p$ in $R$ it must create a $3$-cup or an $l$-cap. We can choose $p$  inside $R$ to have  smaller $x$ coordinate than any element of $P$, hence we can ensure that $p$ is not part of any $3$-cup. Therefore $p$ must be in a $l$-cap and $P$ must have at least $l-1$ elements. 

On the other hand a point set forming an $(l-1)$-cap is semisaturated for $k=3$. Indeed if we add a point, then it either creates a 3-cup or any three points form a 3-cap, which means that the whole set is a cap. The $l=3$ case can be handled similarly.

Now we will assume that $k\ge 5$ and $l\ge 5$; the case when $k=4$ or $l=4$ will be settled later. We can define $L$ and $R$ as above. If we add $p$ anywhere in $R$ it must create a $k$-cup or an $l$-cap. Since $k\ge 4$ and $p$ lies below the lines of $L$, $p$ can only create an $l$-cap and furthermore $p$ is either the first element of this cap or the last one. We can choose $p$ inside $R$ to have a smaller $x$ coordinate than any element of $P$ (see Figure \ref{fig:cupcap}). This ensures that $p$ is the first element in the $l$-cap it has created. Now we move $p$ continuously inside $R$, all the while increasing its $x$-coordinate, until it has bigger $x$-coordinate than any element of $P$. At this point $p$ cannot be the first element of the $l$-cap it creates, thus during this movement there must be a last moment where $p$ is the first element of some $l$-cap containing it. Clearly the change happens as $p$ passes below an element $p_{below}$ of $P$. Let $x_{below}$ denote the $x$-coordinate of this point. 

\begin{figure}[!ht]
    \centering
    \includegraphics[width=10cm]{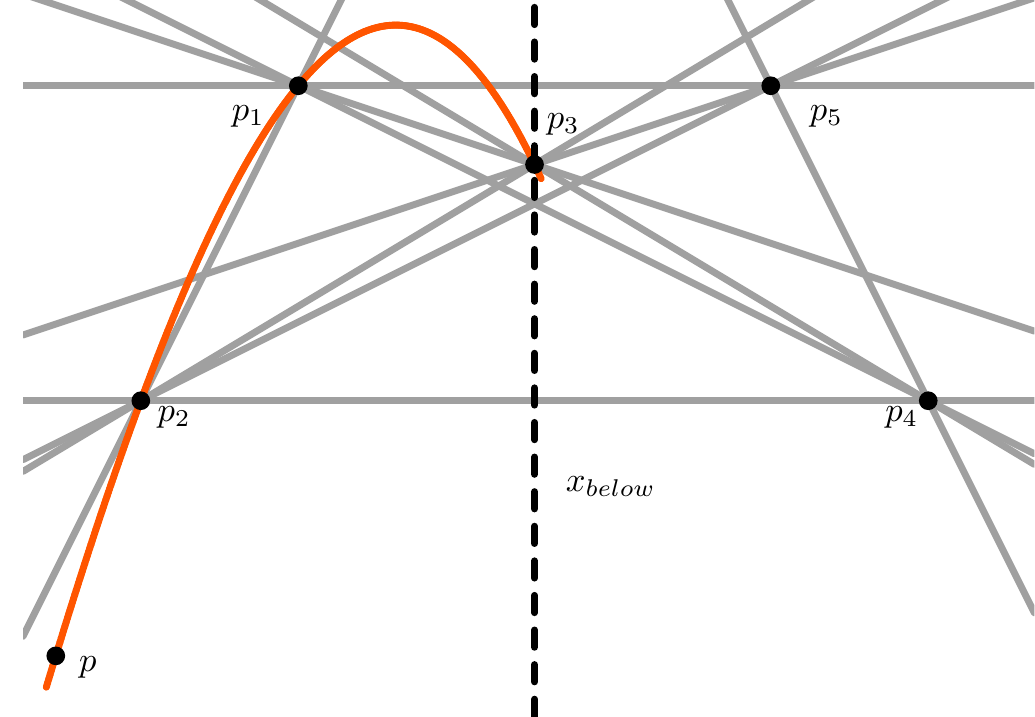}
    \caption{If $p$ is to the right of $x_{below}$, then it is not the first element of any $4$-cap.}
    \label{fig:cupcap}
\end{figure}

If we put $p$ in $R$ slightly after $x_{below}$, then it must extend some $(l-1)$-cap $A_1$ whose points lie to the left of $x_{below}$, except maybe for $p_{below}$. Similarly if we put it slightly before $x_{below}$, then it must extend some $(l-1)$-cap $A_2$ that lies to the right of $x_{below}$. 

In the same way we can define $p_{above}$, $x_{above}$ and two $(k-1)$-cups $U_1$ and $U_2$ such that they lie on the left and right side of $x_{above}$. 

In total we have found $|A_1\cup A_2 \cup U_1 \cup U_2|$ elements. Clearly $|A_1\cap A_2|\le 1$ and $|U_1\cap U_2|\le 1$.  Since any cup intersects any cap in at most two points we have $|A_i\cap U_j|\le 2$. Also either $x_{below}< x_{above}$ or $x_{below}> x_{above}$ or $x_{below}= x_{above}$. In the first case $A_1\cap U_2=\emptyset$ and in the second case $A_2\cap U_1=\emptyset$, giving us 
\[|A_1\cup A_2 \cup U_1 \cup U_2|\ge 2k-2+2l-2-|A_1\cap A_2|-|U_1\cap U_2|-3\cdot 2\ge 2k+2l-12.\]
If $x_{below}= x_{above}$ we have $|A_1\cap U_2|\le 1$ and $|A_2\cap U_1|\le 1$  so we have 
\[|A_1\cup A_2 \cup U_1 \cup U_2|\ge 2k-2+2l-2-|A_1\cap A_2|-|U_1\cap U_2|-2\cdot 2-1-1\ge 2k+2l-12.\] 

For the upper bound we give a construction. First consider the point set shown in Figure~\ref{fig:cupcapconst}.

\begin{figure}[!ht]
    \centering
    \includegraphics[width=7cm]{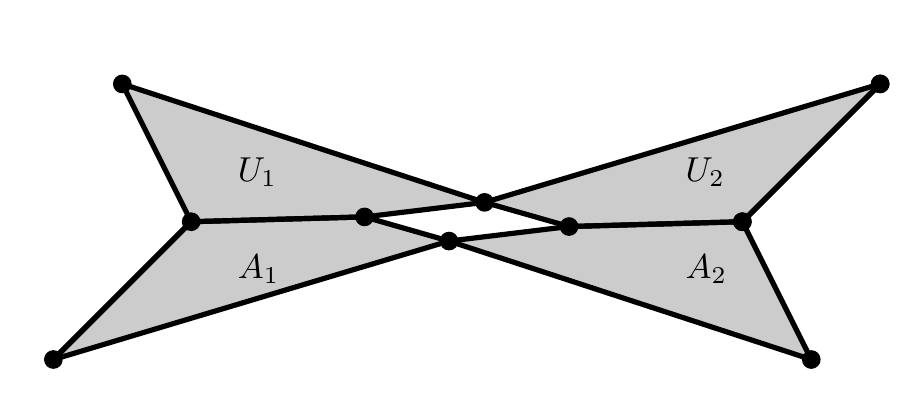}
    \includegraphics[width=7cm]{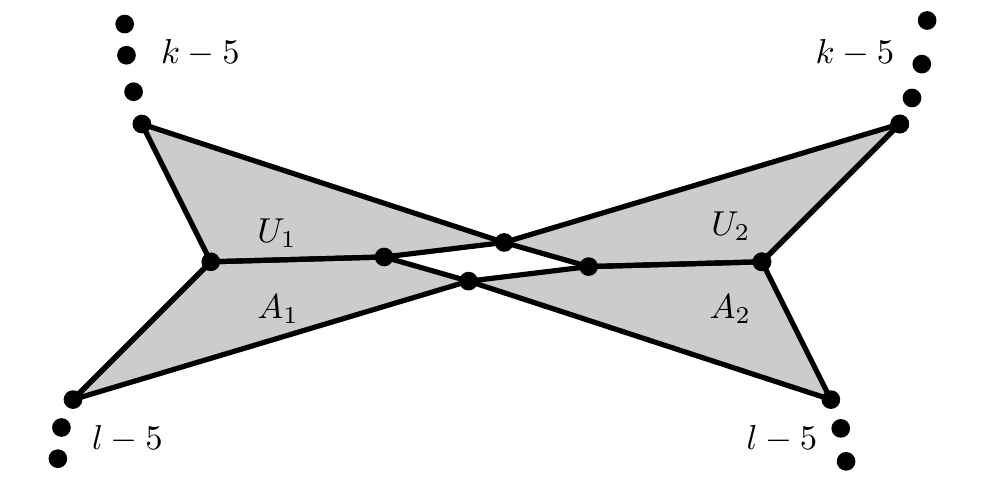}
    \caption{ Cup-cap semisaturation for $k=l=5$ and for $k=8,l=7$}
    \label{fig:cupcapconst}
\end{figure}

This point set consists of 10 points, and it is saturated for 5-cups and 5-caps. We show this by dividing the plane into regions and for each region showing four points of the point set that form either a 5-cup or a 5-cap with any point of the region.  Consider the regions in Figure~\ref{fig:cupcapregions}.

\begin{figure}[!ht]
    \centering
    \includegraphics[width=0.3\textwidth]{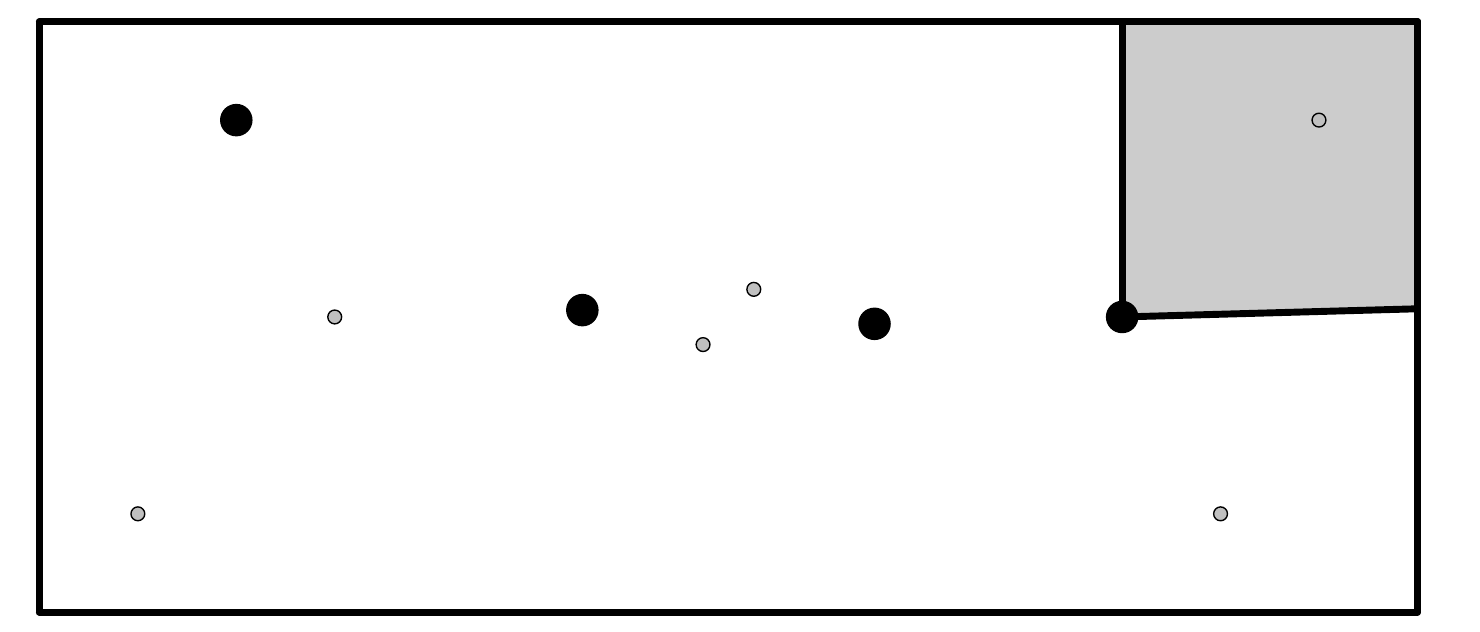}
    \includegraphics[width=0.3\textwidth]{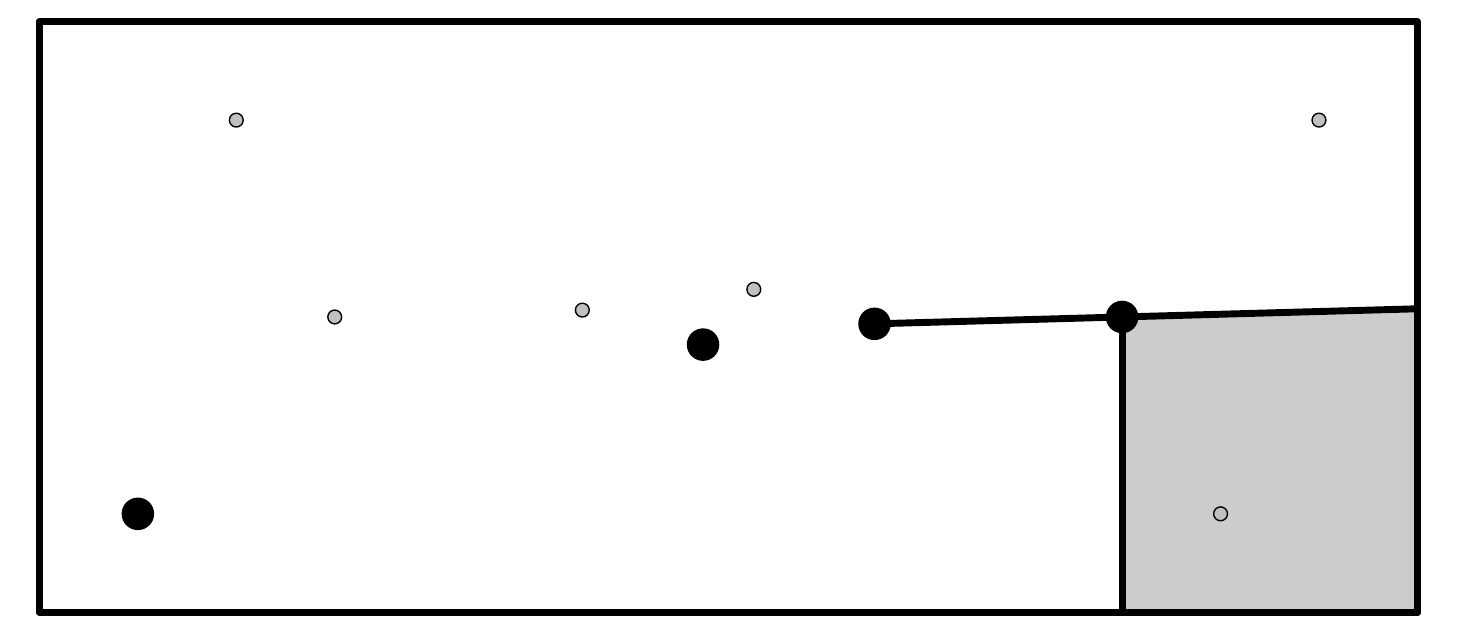}
    \includegraphics[width=0.3\textwidth]{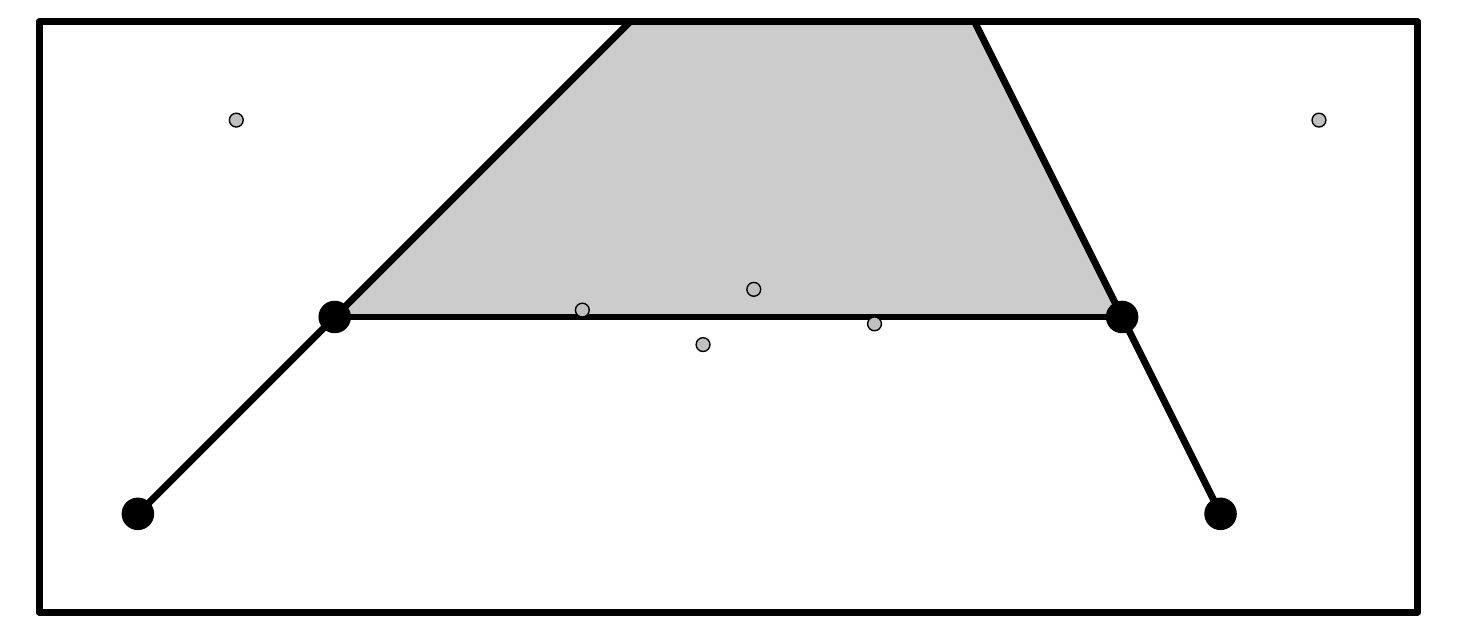}
    \includegraphics[width=0.3\textwidth]{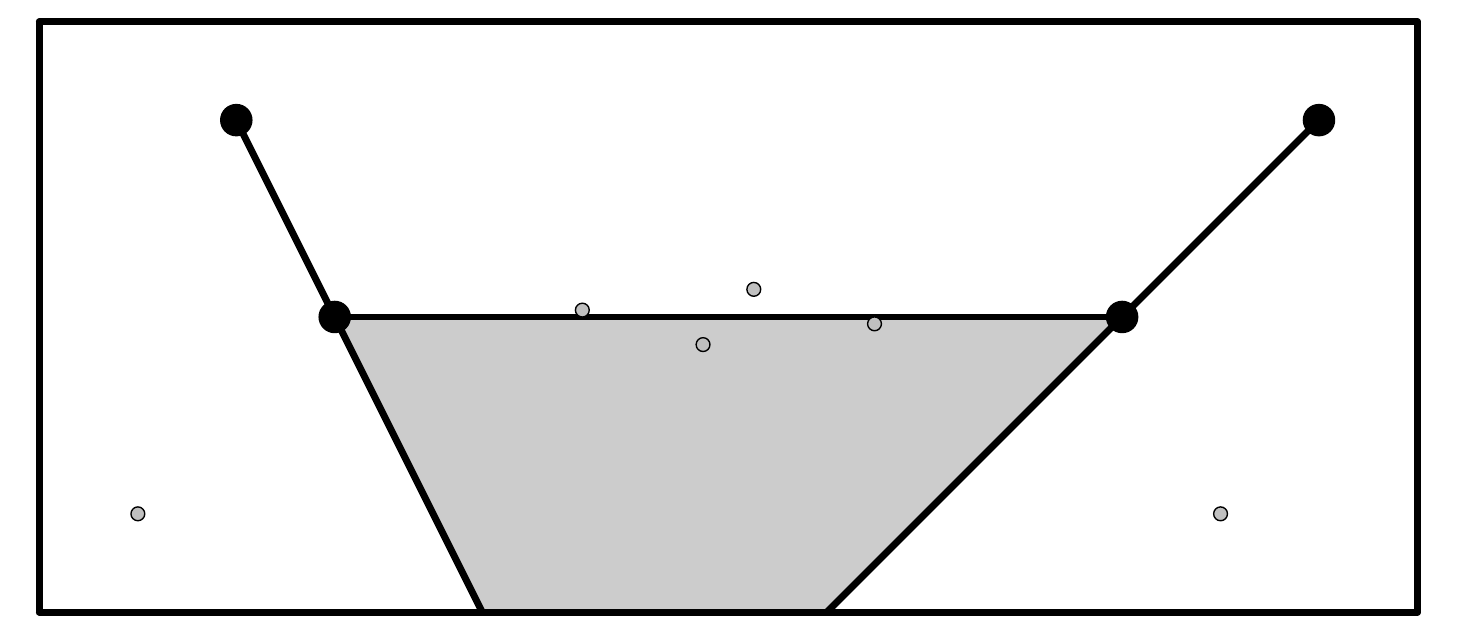}
    \includegraphics[width=0.3\textwidth]{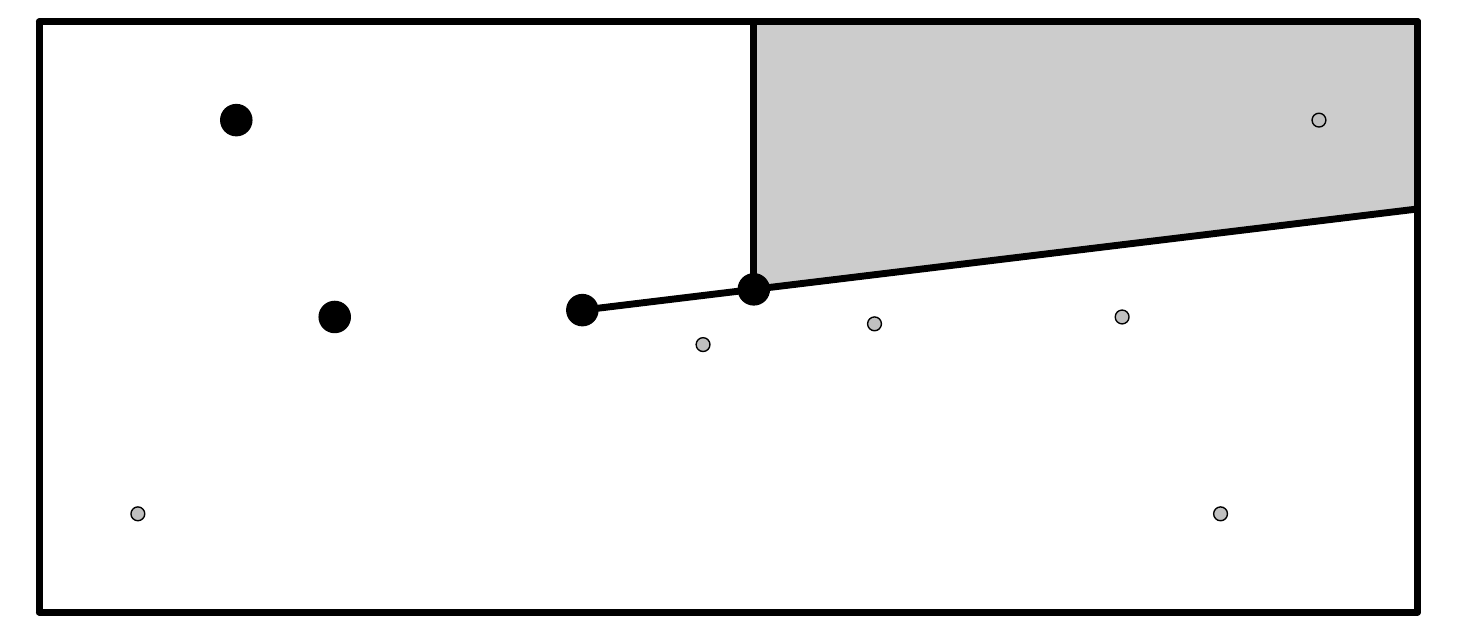}
    \includegraphics[width=0.3\textwidth]{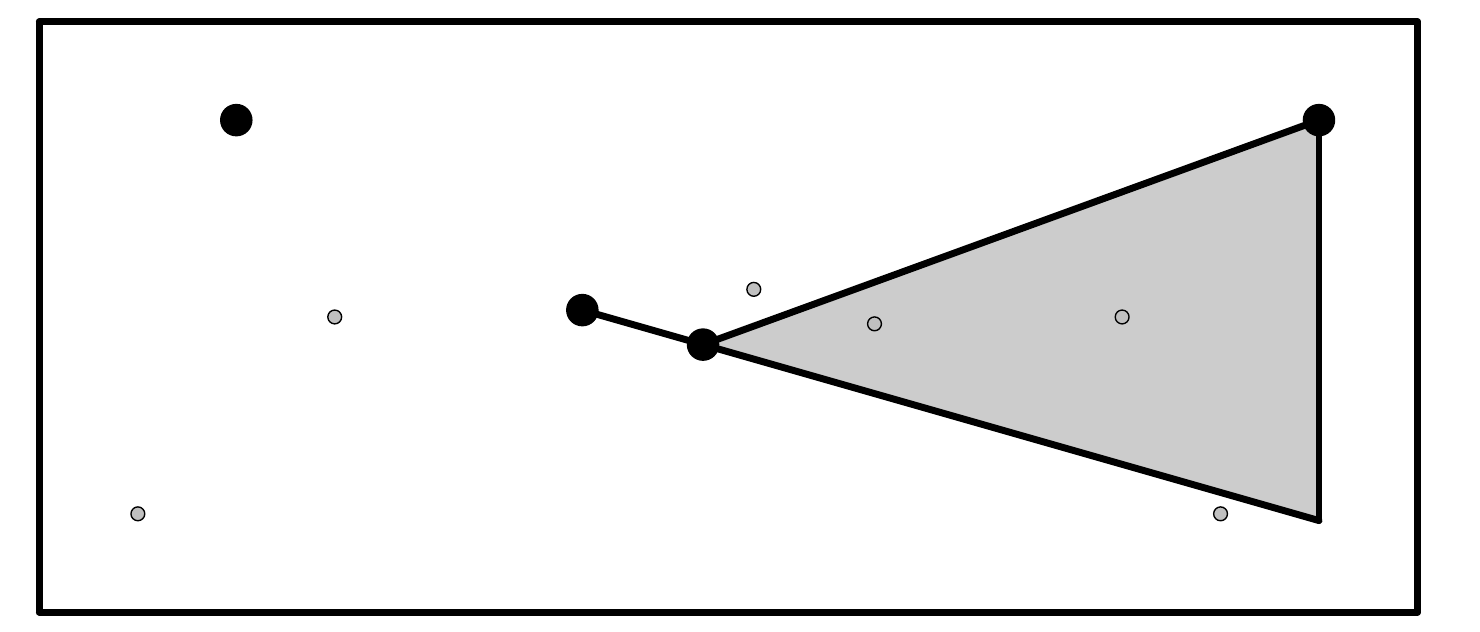}
    \includegraphics[width=0.3\textwidth]{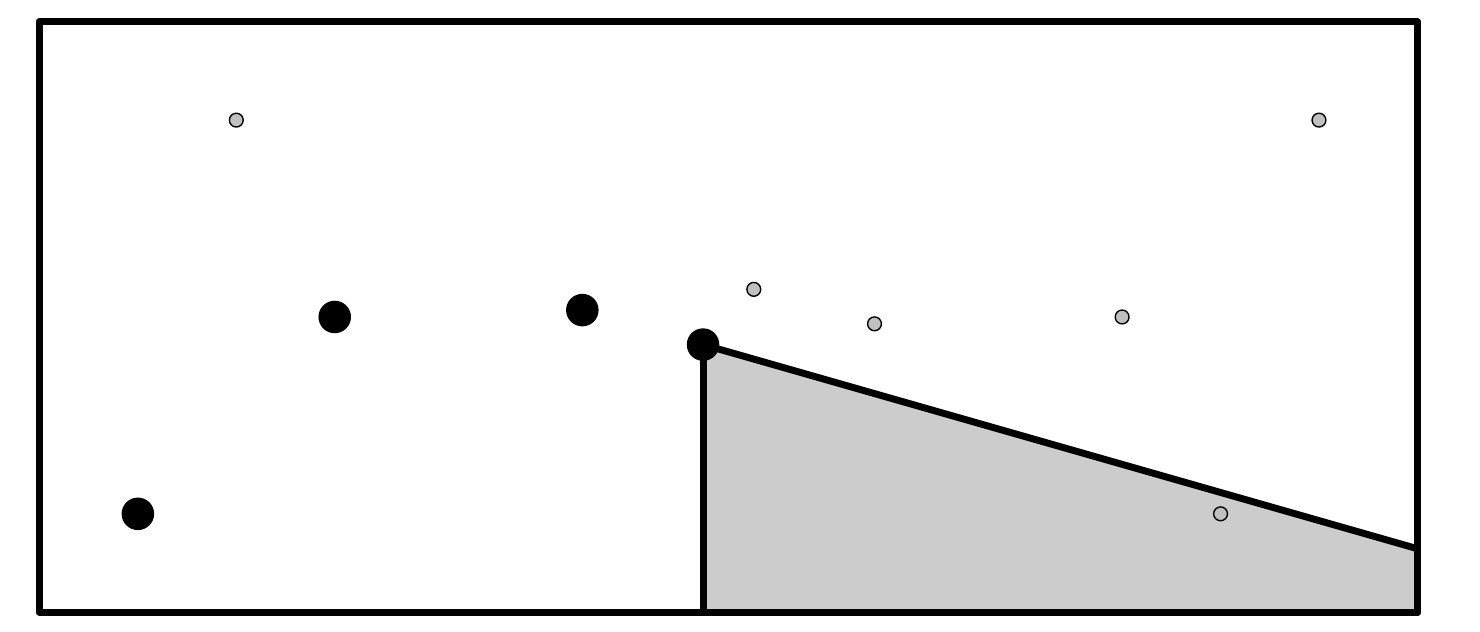}
    \includegraphics[width=0.3\textwidth]{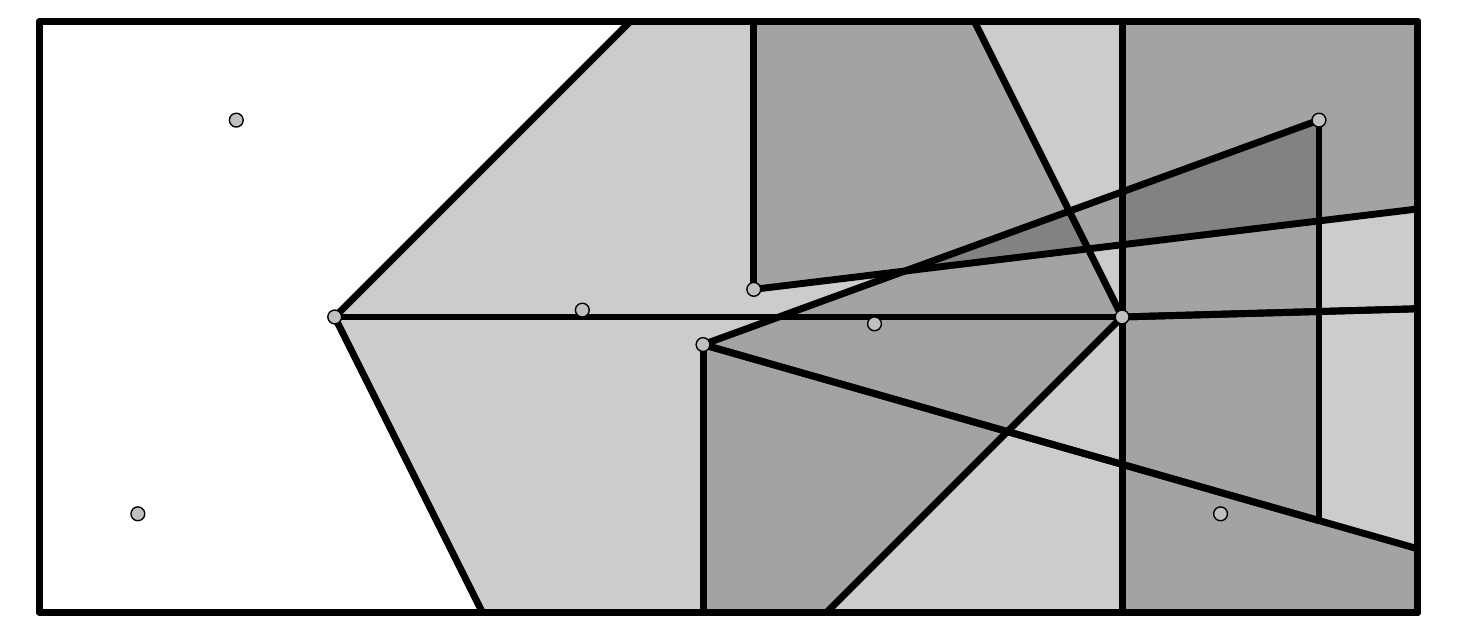}
    \caption{ The regions blocked by the 4-cups and 4-caps.}
    \label{fig:cupcapregions}
\end{figure}

In the first seven subfigures we have drawn a region and indicated which four points of the point set blocks that region. In the eight subfigure we have drawn all the regions. As we can see these regions cover half of the plane. Since the point set is centrally symmetric, this is enough. Hence for $k=l=5$ we have semisaturation with $2k+2l-10$ points. 

Now we will show that this construction can be extended for $k,l\ge 5$. In Figure \ref{fig:cupcapconst} we can see $A_1$, $A_2$, $U_1$ and $U_2$. To get a construction for $k,l$ we just extend $A_1$ to the left and $A_2$ to the right with $l-5$ elements and $U_1$ to the left and $U_2$ to the right with $k-5$ elements. See Figure \ref{fig:cupcapconst} for an example. The resulting configuration will be semisaturated. Considering the same regions as in Figure \ref{fig:cupcapregions} will work. If a region was blocked by a $4$-cup, that $4$-cup is now extended by $k-5$ elements, so we have a blocking $(k-1)$-cup. Similarly if a region was blocked by a $4$-cap, that $4$-cap is now extended by $k-5$ elements, so we have a blocking $(k-1)$-cap. Therefore we have found a semisaturated set with $10+2(k-5)+2(l-5)=2k+2l-10$ points. 

In the case of $l=4$ and $k\ge 5$ the construction is quite similar. A possible configuration for the  $k=5$ case is given in Figure \ref{fig:4case} and the blocked regions are given in Figure \ref{fig:4caseall}. For $k>5$ we can extend $U_1$ and $U_2$ just as we did in Figure \ref{fig:cupcapconst}. We leave the details for the interested readers.  

\begin{figure}[!ht]
    \centering
    \includegraphics[width=0.3\textwidth]{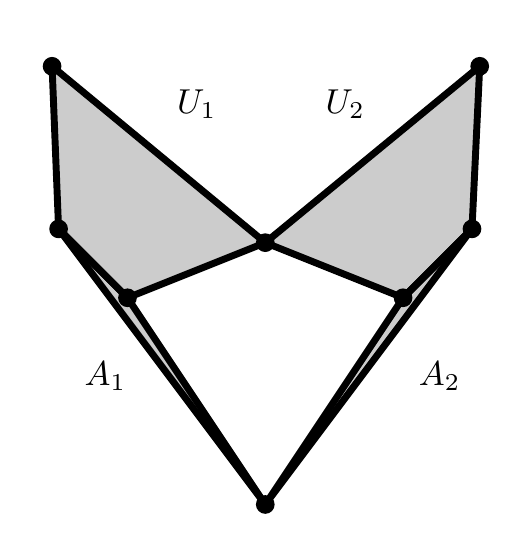}
    \includegraphics[width=0.39\textwidth]{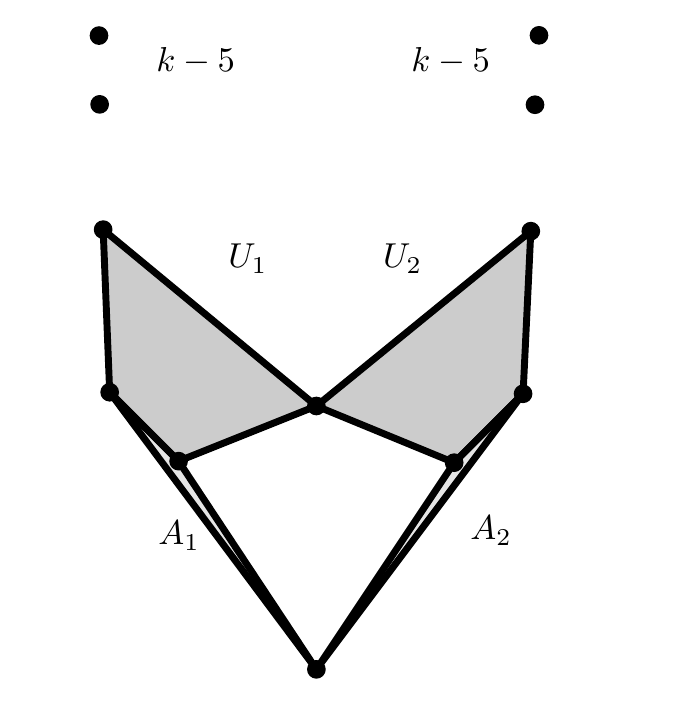}
    \caption{Cup-cap semisaturation for $l=4, k = 5$ and for $l=4, k=7$.}
    \label{fig:4case}
\end{figure}

\begin{figure}[!ht]
    \centering
    \includegraphics[width=0.15\textwidth]{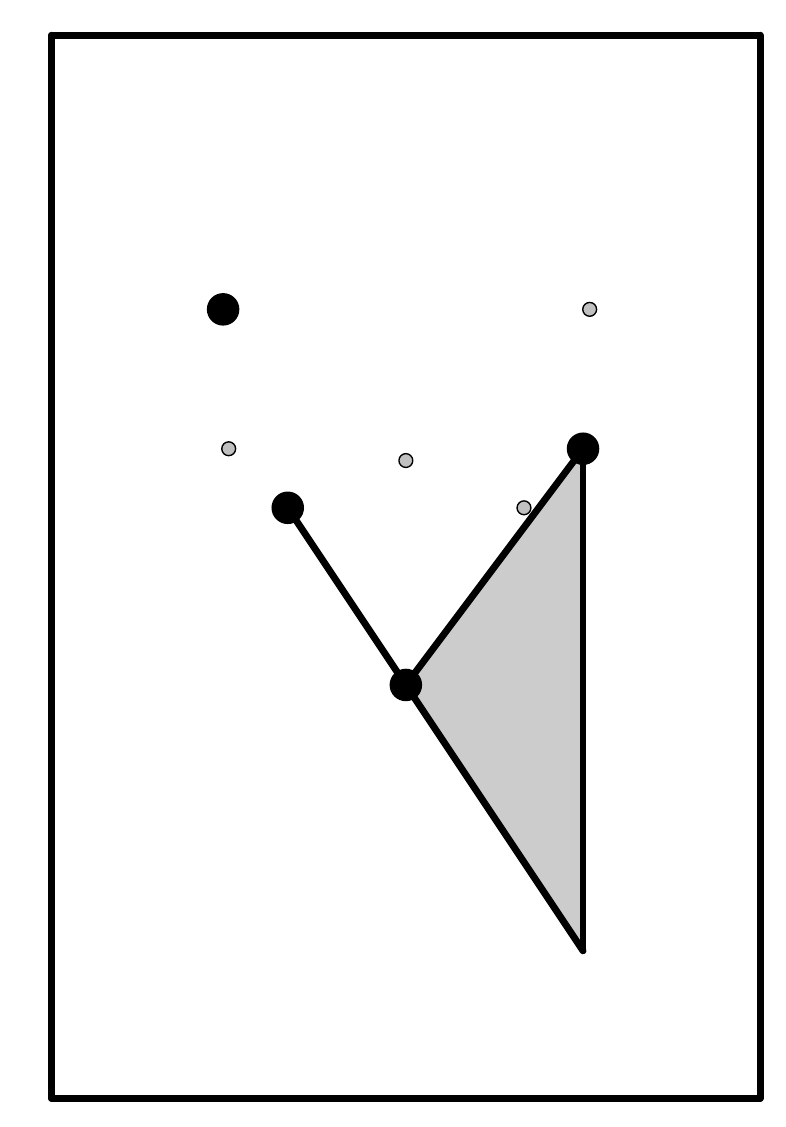}
    \includegraphics[width=0.15\textwidth]{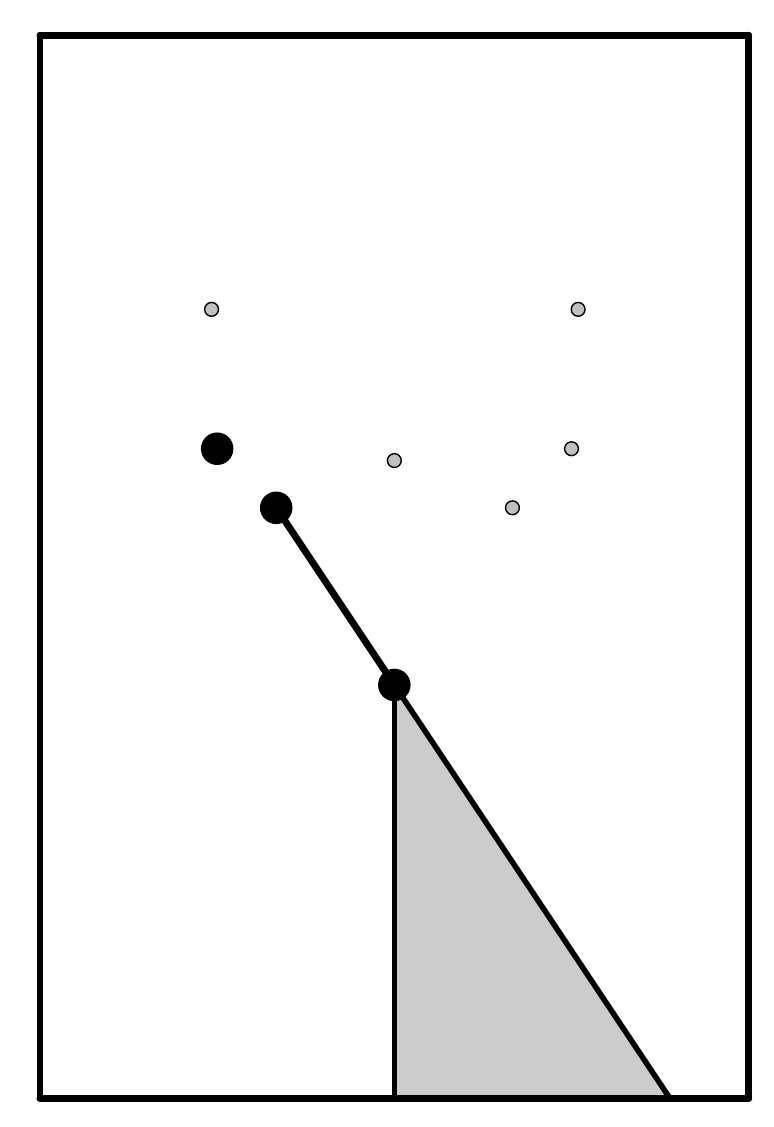}
    \includegraphics[width=0.15\textwidth]{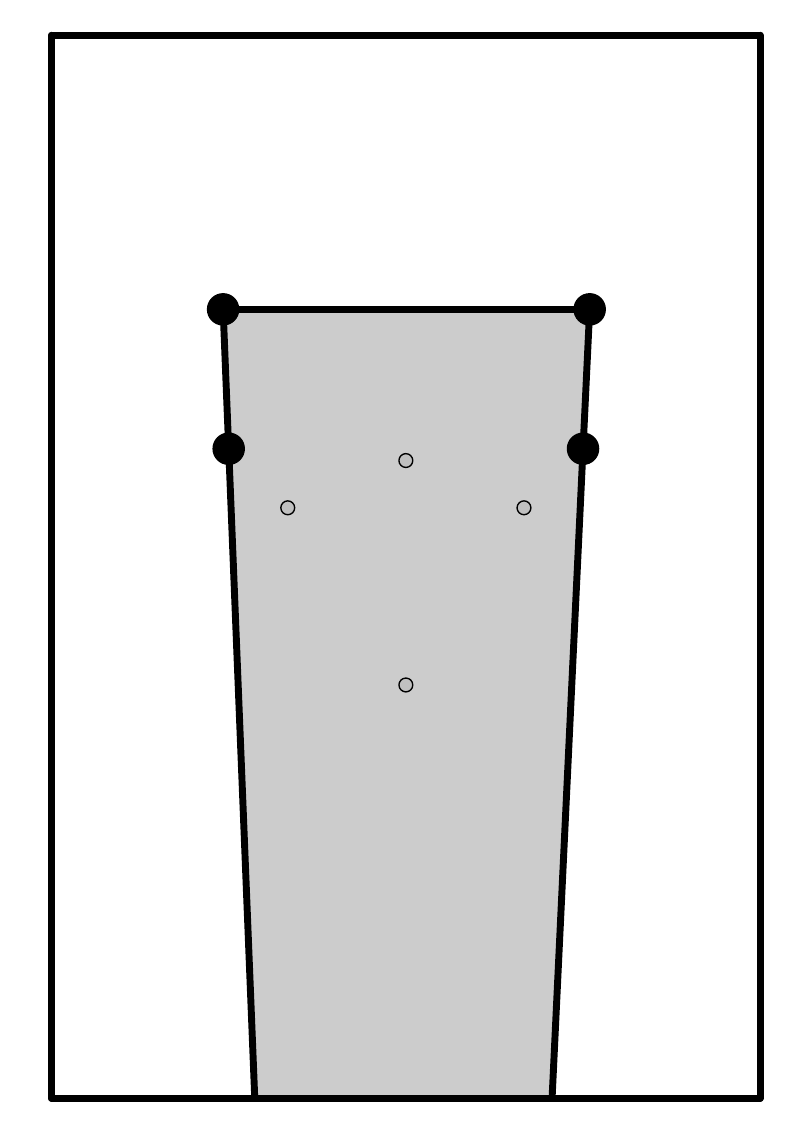}
    \includegraphics[width=0.15\textwidth]{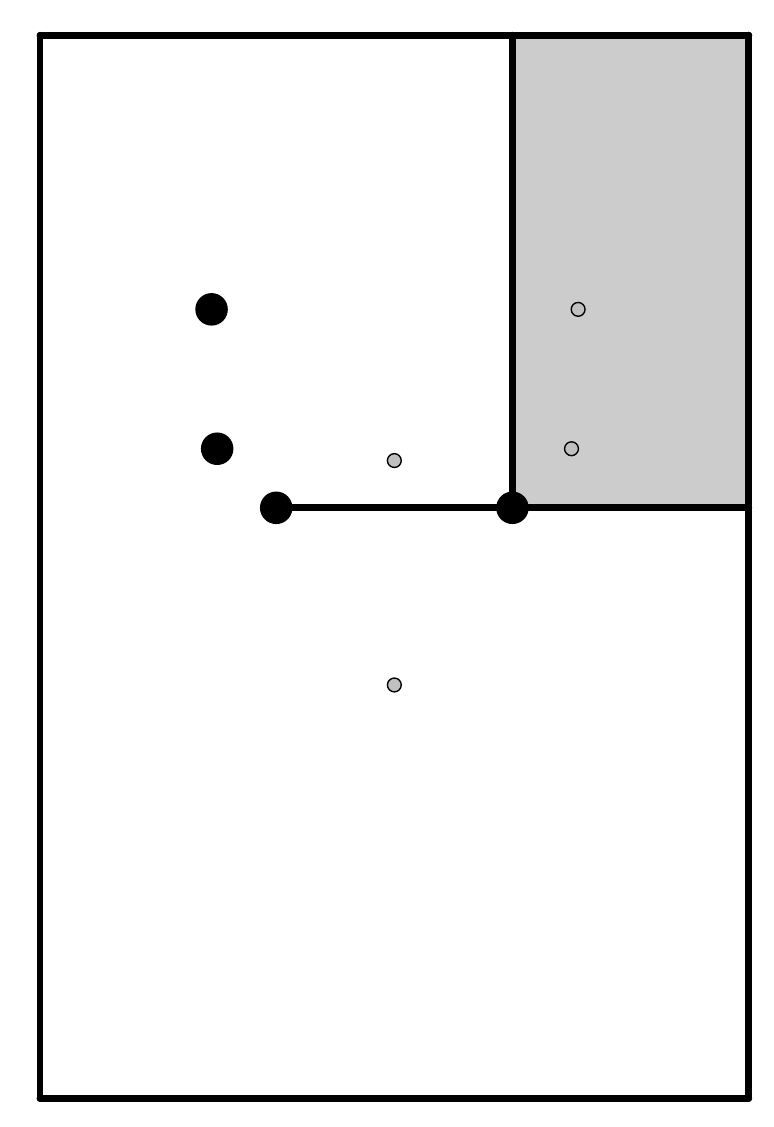}
    \includegraphics[width=0.15\textwidth]{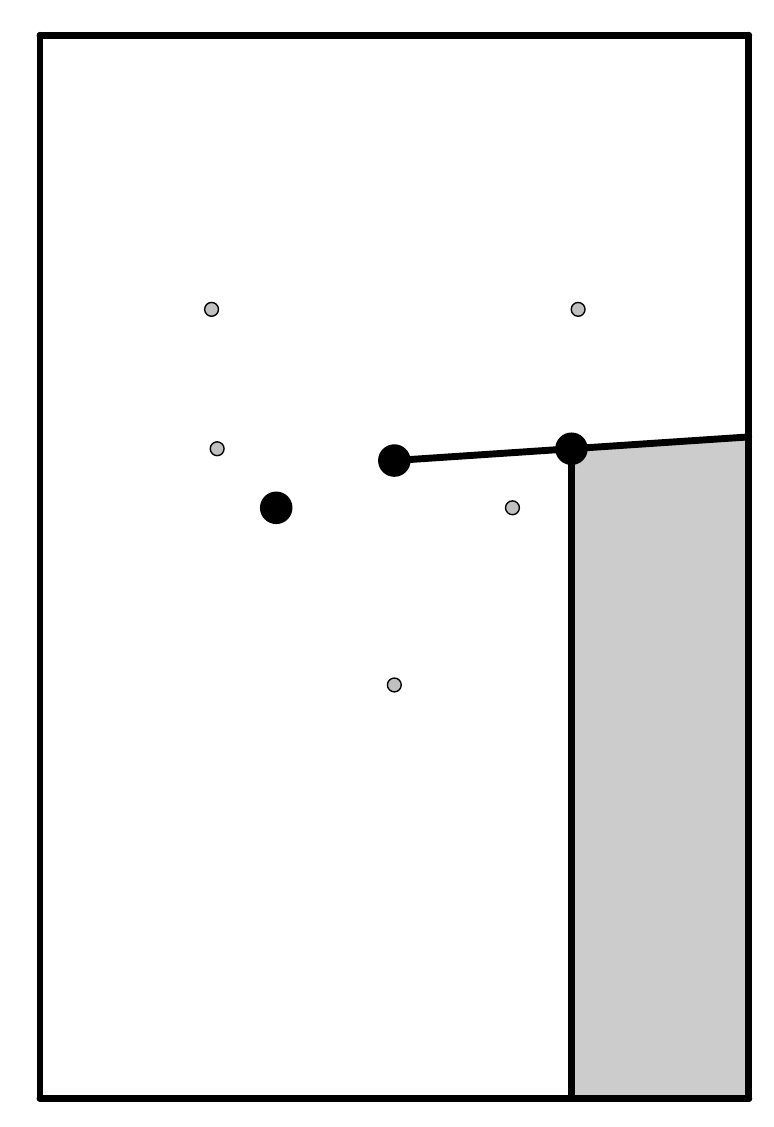}
    \includegraphics[width=0.15\textwidth]{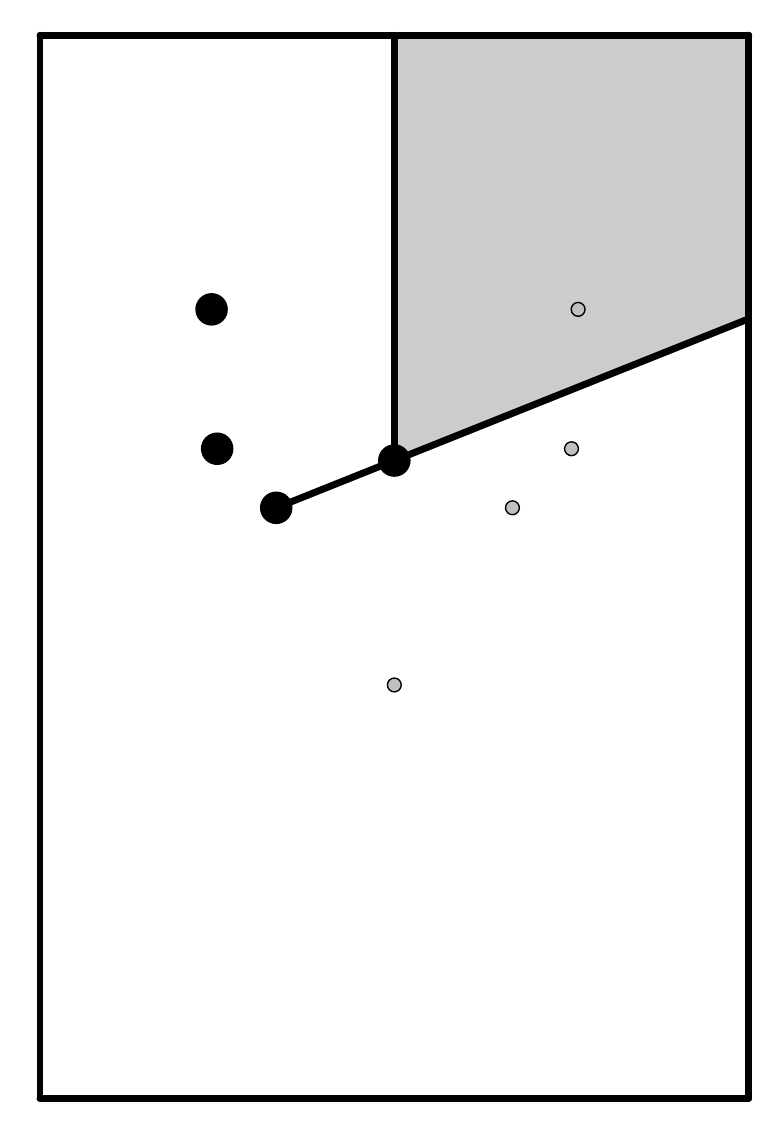}
    \includegraphics[width=0.15\textwidth]{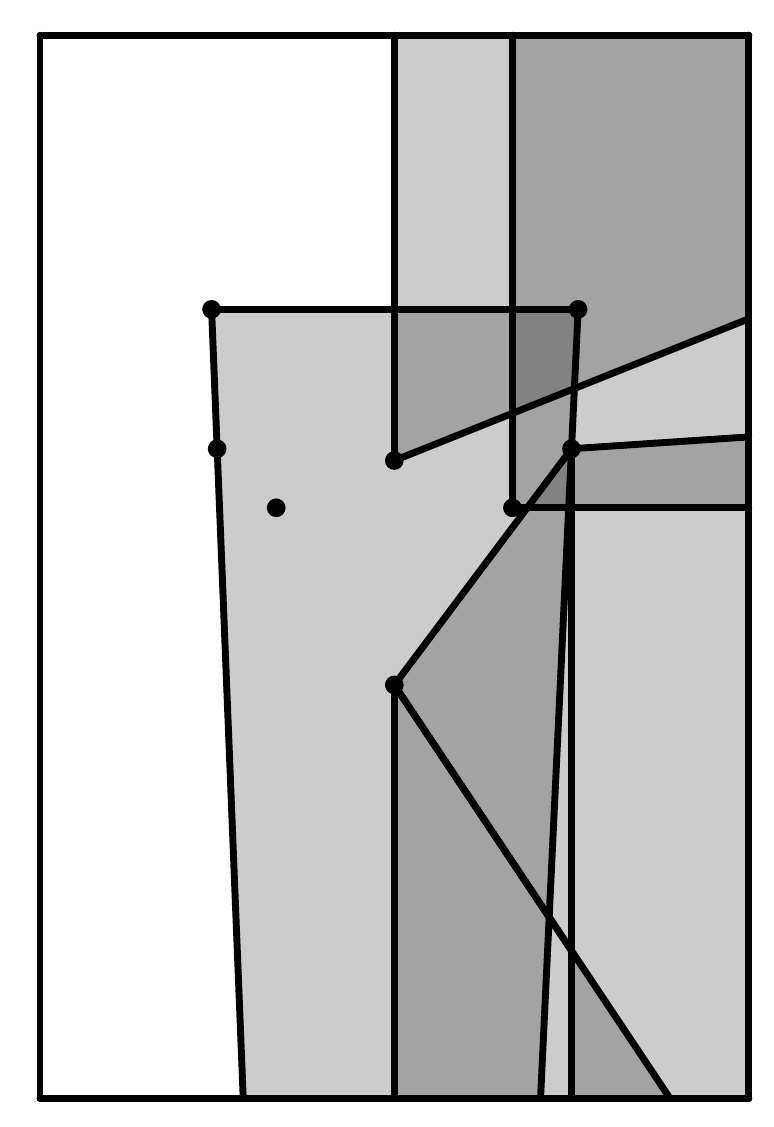}
    \caption{The regions blocked by the 4-cups and 3-caps.}
    \label{fig:4caseall}
\end{figure}

Next we  show that at least $2k-2$ points are required to obtain a saturated construction. Indeed, by the same reasoning as in the $l,k \ge 5$ case we can find two cups $U_1,U_2$ of size $k-1$ intersecting in at most one point and two caps $A_1,A_2$ of size $3$ intersecting in at most one point. If $U_1$ and $U_2$ are disjoint, then we have already found $2k-2$ many points, so suppose they intersect in one point $q$. We know that either $A_1$ lies to the left of $q$ or $A_2$ lies to the right of $q$ (it is possible that $A_1$ or $A_2$ contains $q$). In either case we must have one more point since no three points of $U_1$ nor of $U_2$ form a cap.  
\end{proof}

Now we continue with the semisaturation of convex point sets.	

\begin{proof}[Proof of Theorem~\ref{thm:osat00}]
Suppose that to the contrary there is a semisaturated set of points $S$ with $n-1 + \floor{\frac{n-2}{d}}-s$ points in $\R^d$, for some $s\ge 1$.  Denote by $S_1,S_2,\dots,S_m$ the subsets of $S$ that are convex $(n-1)$-sets.  If $\cap_{i=1}^m \interior(\conv(S_i)) \neq \varnothing$, then we may add any point in the intersection without yielding $n$ points in convex position. We can also add this point such that the resulting point set is in general position.  The interior of a convex set is convex, hence by Helly's theorem it is sufficient to show that the intersection of any $(d+1)$ of the sets from  $\{\interior(\conv(S_i))\}_{i\in [m]}$ is nonempty.  

Consider $d+1$ sets in $\{S_i\}_{i\in [m]}$: $S_{i_1},S_{i_2},\dots,S_{i_{d+1}}$.  They each have size $n-1$ and are contained in the point set $S$ of size $n-1 + \floor{\frac{n-2}{d}}-s$, thus we have that 

\begin{align*}
\abs{(S_{i_1}\cap S_{i_2}\cap \dots \cap S_{i_{d+1}})^c} &= \abs{S_{i_1}^c\cup S_{i_2}^c\cup \dots \cup S_{i_{d+1}}^c}  \\& \leq \abs{S_{i_1}^c}+ \abs{S_{i_2}^c} + \dots + \abs{S_{i_{d+1}}^c}  =(d+1)\left(\floor{\frac{n-2}{d}}-s\right).
\end{align*}
Therefore 
 \begin{align*}
 \abs{S_{i_1}\cap S_{i_2}\cap \dots \cap S_{i_{d+1}}}& \ge \left(n-1+\floor{\frac{n-2}{d}}-s\right)-(d+1)\left(\floor{\frac{n-2}{d}}-s\right) \\&= n-1-d\left(\floor{\frac{n-2}{d}}-s\right)\ge1+ds\ge d+1.
 \end{align*}

Since the original point set is in general position, these $d+1$ points in the intersection span a non-degenerate simplex. It follows that the interiors of $\conv(S_{i_1}),\dots, \conv(S_{i_{d+1}})$ intersect, as required. 
\end{proof}

\begin{proof}[Proof of Theorem~\ref{thm:convex}]
We will construct a set of points $S$ such that $S$ is semisaturated in the plane $\mathbb{R}^2$. Consider a convex polygon $Q =\conv(v_0,v_1,\ldots, v_{2n-3})$ with $2n-4$ sides such that the side $v_i v_{i+1}$ is parallel to the side $v_{n+i-2} v_{n+i-1}$ (where the indices are modulo $2n-4$). For ease of reference, we call the pair of sides $v_i v_{i+1}$, $v_{n+i-2} v_{n+i-1}$ \textit{opposite sides} of $Q$. Let $S=\{v_0,v_1,\dots,v_{2n-3}\}$ be the vertex set of $Q$. We claim that $S$ is semisaturated in $\mathbb{R}^2$.

\begin{claim}\label{cl:inside}
Let $P$ be any point contained in $Q$ such that $Q\cup\{P\}$ is in general position. Then $S \cup \{P\}$ has a convex $n$-gon with $P$ as one of its vertices.
\begin{proof}
Let $P$ be an arbitrary point in the interior of $Q$. Consider the chord $v_0v_{n-2}$ (see Figure~\ref{fig:convexfig}). This divides $Q$ into two $(n-1)$-gons $\{ v_{0}, v_{1},\ldots, v_{n-2}, \}$ and $\{ v_{n-2}, v_{n-1},\ldots, v_{0}, \}$. Since the points are in general position, $P$ must lie in the interior of one of these $(n-1)$-gons. Then $P$ and the points of the other $(n-1)$-gon form a convex $n$-set.
\end{proof}
\end{claim}

\begin{claim}\label{cl:outside}
Let $P$ be any point not contained in $Q$ such that $Q\cup\{P\}$ is in general position. Then $S \cup \{P\}$ has a convex $n$-gon with $P$ as one of its vertices.
\end{claim}

We say that a point $A\in Q$ can be \emph{seen} from $P$ if $\overline{PA}\cap Q={A}$. A side $v_i,v_{i+1}$ can be seen from $P$ if all of its points can be seen from $P$. Clearly $P$ can see at most one from a pair of opposite sides. So there are at least $n-2$ sides that cannot be seen from $P$. Each of these sides will be a side of the convex hull $\conv(Q\cup\{P\})$. There will be two additional sides of this convex~hull incident to $P$, so $\conv(Q\cup\{P\})$ has at least $n$ sides. Therefore we can choose $n-1$ points from $S$ such that they form a convex $n$-gon with $P$.   

\begin{figure}[!ht]
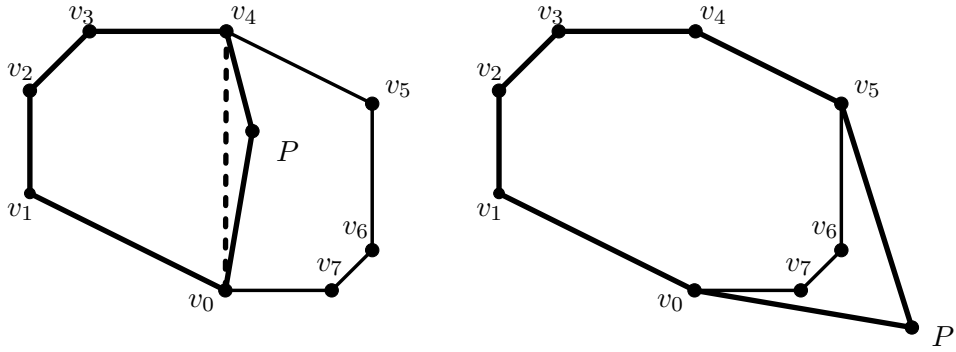

    \centering
    \include{convex}
    \caption{Finding a convex $n$-gon in the extended point set.}
    \label{fig:convexfig}
\end{figure}

Claim~\ref{cl:inside} and Claim~\ref{cl:outside} imply that $S$ is semisaturated, hence $\osat_{\C}(n) \leq 2n-4$.

Now we prove the lower bound. 
Suppose $S$ is a semisaturated set of points in the plane. 
The set $S$ determines $\binom{|S|}{2}$ lines, which partition the plane into regions. Let $P_1$ be a point in one of the infinite regions (and not between any pair of parallel lines), and let $P_2$ be another point in the opposite infinite region. That is $P_1$ and $P_2$ lie on different sides for each of the $\binom{|S|}{2}$ lines.

Since $S$ is semisaturated, there are two sets $S_1,S_2\subset S$ such that $S_1\cup \{P_1\}$ and $S_2\cup \{P_2\}$ are convex $n$-gons. We claim that $|S_1\cap S_2|\le 2$. 
Suppose $v_1,v_2,v_3 \in S_1\cap S_2$, and assume that $v_1,v_2,v_3,P_1$ is a convex quadrilateral (in that order). Observe that $P_1$ is in the same side as $v_3$ with respect to the line $v_1v_2$ and on the same side as $v_1$ with respect the line $v_2v_3$. On the other hand, $P_2$ is in the opposite side of $v_3$ with respect the line $v_1v_2$ and on the opposite of $v_1$ with respect the line $v_2v_3$.  The points $v_1,v_2,v_3,P_2$ cannot form a convex quadrilateral (in any order). Indeed, one of the sides of this quadrilateral must be either $v_1v_2$ or $v_2v_3$, but the line defined by $v_1v_2$ would separate $P_2$ from $v_3$, and the line defined by $v_2v_3$ would separate $P_2$ from $v_1$. This cannot happen since $S_2 \cup \{P_2\}$ is in convex position, a contradiction.

\begin{figure}[!ht]
    \centering
    \includegraphics{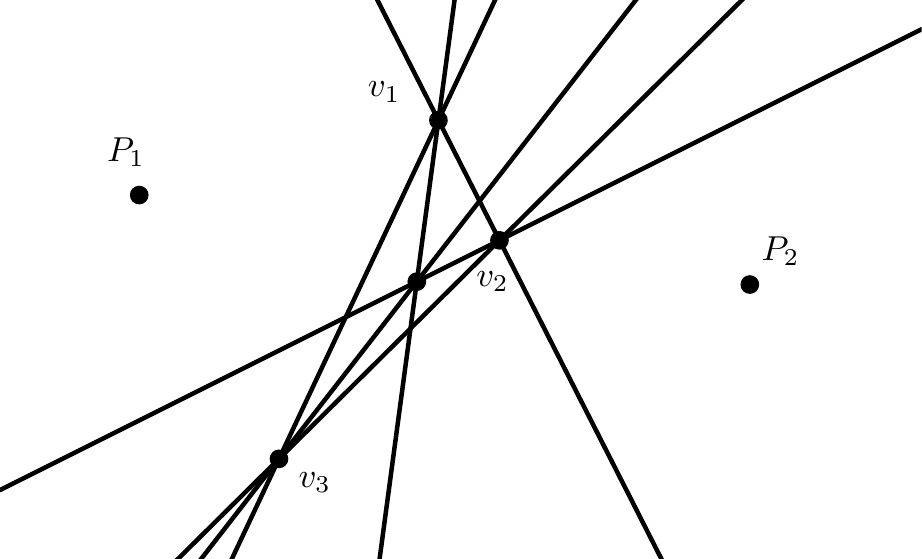}
    \caption{Convex $n$-gons containing $P_1$ and $P_2$ intersect in at most two points.}
    \label{fig:conv2}
\end{figure}
Hence $|S_1\cup S_2|=|S_1|+|S_2|-|S_1\cap S_2|\ge n-1+n-1-2=2n-4.$
\end{proof}

\section{General treatment of saturation questions}\label{general}

In this section we provide a general formulation for many of the the problems we have considered in this paper. Given a $c$-edge-colored complete $s$-uniform hypergraph $H=(V,E)$, we say that a subset of vertices $S \subseteq V$ forms a monochromatic complete subhypergraph of $H$ in color $i$ if the ($s$-uniform) subhypergraph induced by $S$ has only hyperedges in color $i$.

Many of the problems we considered have the following form.
\begin{defi}\label{gengen}
     	Given constants $c$ and $s$, let $\F_0$ be the family of complete $s$-uniform hypergraphs whose edges are colored with $c$ colors (numbered by $1,2,\dots, c$). For a subfamily $\F$ of $\F_0$, a member $F$ of $\F$ is \emph{saturated} if for every $i$, $F$ does not contain a monochromatic complete subhypergraph of size $k_i$ and color $i$, but every $F' \in \F$ that extends $F$ contains a monochromatic complete subhypergraph of size $k_i$ of color $i$ for some $i$. $F$ is \emph{semisaturated} if we omit the first condition,  that is, if $F \in\F$ and every $F'\in \F$ that extends $F$ contains a monochromatic complete subhypergraph of size $k_i$ of color $i$ for some $i$ which is not in $F$.
     	
     	Let $\ram_{\F}(k_1,\dots k_c)$ denote the size (number of vertices) of the largest saturated $F\in \F$, and let $\sat_{\F}(k_1,\dots k_c)$ denote the size of the smallest saturated $F\in \F$. Finally, let $\osat_{\F}(k_1,\dots k_c)$ denote the size of the smallest semisaturated $F\in \F$.	     	
\end{defi}

\begin{obs}
	For any $\F$ and positive integers $k_1, \ldots, k_c$, 
	\[\osat_{\F}(k_1, \ldots, k_c)\le \sat_{\F}(k_1, \ldots, k_c) \le \ram_{\F}(,k_1, \ldots, k_c).\]
\end{obs}

Note that whenever $\sat=\ram$ holds, all saturated members of $\F$ have the same size. Thus we gain further insight into the respective Ramsey-type problem as well.
Moreover, when $c=2$, one can regard the problem as the first color class forming an (uncolored) hypergraph $H$. Then it follows that a complete subhypergraph in the first color is a complete subhypergraph in $H$ while a complete subhypergraph in the second color is an independent set in $H$.

Definition~\ref{gengen} is quite general. In this paper we have introduced saturation problems for graphs, posets, monotone point sets, and cups and caps. All of these fit into this formulation. First, the graph case we get by setting $s=2$ and $\F=\F_0$. We get the poset case by setting $c=s=2$ and letting $\F$ be the family of those $2$-edge-colored graphs that we can obtain as the comparability graph of a poset. We obtain the monotone point set case by setting $c=s=2$ and letting $\F$ be the family of those $2$-edge-colored graphs that we can obtain from the pairs of elements in a sequence by coloring the increasing pairs red and the decreasing pairs blue. Finally, the cup and cap case we get by setting $c=2$, $s=3$ and letting $\F$ be the family of those $2$-colored complete $3$-uniform hypergraphs that we can obtain by taking a point set in general position and coloring a  triple red if it forms a cup and blue if it forms a cap (note that every triple forms a cup or a cap). The only problem we considered that does not fit into this formulation is the case of convex subsets of points.

It is interesting that for both the $2$-colored graph case and the poset case we have $\sat(k,l)=(k-1)(l-1)$, yet we could not find any general reasoning that handles both of these cases at once. It is interesting that the relative behavior of $\sat$, $\osat$ and $\ram$ can vary substantially depending on the setting. Indeed, for graphs $\osat=\sat$ yet $\ram$ is exponential, while for posets and monotone point sets, $\sat$ equals $\ram$ yet $\osat$ is smaller ($\osat$ behaves differently in the latter two cases).

\section{Acknowledgements}
We thank Tuan Tran for pointing out an inaccuracy in Lemma~\ref{decomp} in an earlier version of this manuscript.

 \end{document}

%% file: posetconst4.tex
\begin{tikzpicture}[line cap=round,line join=round,>=triangle 45,x=1.0cm,y=1.0cm]
\clip(1.4992778219438057,5.225021719891174) rectangle (5.9350014719574355,9.808057225440328);
\draw [shift={(2.5,7.5)},line width=1.3pt]  plot[domain=1.5707963267948966:3.141592653589793,variable=\t]({1.*0.3*cos(\t r)+0.*0.3*sin(\t r)},{0.*0.3*cos(\t r)+1.*0.3*sin(\t r)});
\draw [shift={(2.5,7.5)},line width=1.3pt]  plot[domain=3.141592653589793:4.71238898038469,variable=\t]({1.*0.3*cos(\t r)+0.*0.3*sin(\t r)},{0.*0.3*cos(\t r)+1.*0.3*sin(\t r)});
\draw [shift={(5.5,7.5)},line width=1.3pt]  plot[domain=-1.5707963267948966:0.,variable=\t]({1.*0.3*cos(\t r)+0.*0.3*sin(\t r)},{0.*0.3*cos(\t r)+1.*0.3*sin(\t r)});
\draw [shift={(5.5,7.5)},line width=1.3pt]  plot[domain=0.:1.5707963267948966,variable=\t]({1.*0.3*cos(\t r)+0.*0.3*sin(\t r)},{0.*0.3*cos(\t r)+1.*0.3*sin(\t r)});
\draw [line width=1.3pt] (2.2,7.5)-- (2.2,7.5);
\draw [line width=1.3pt] (2.5,7.8)-- (5.5,7.8);
\draw [line width=1.3pt] (5.8,7.5)-- (5.8,7.5);
\draw [line width=1.3pt] (5.5,7.2)-- (2.5,7.2);
\draw [line width=1.3pt] (4.,5.5)-- (4.,6.);
\draw [line width=1.3pt] (4.,6.)-- (4.,6.5);
\draw [line width=1.3pt] (4.,6.5)-- (4.,7.);
\draw [line width=1.3pt] (4.,7.)-- (2.5,7.5);
\draw [line width=1.3pt] (4.,7.)-- (3.5,7.5);
\draw [line width=1.3pt] (4.,7.)-- (4.5,7.5);
\draw [line width=1.3pt] (4.,7.)-- (5.5,7.5);
\draw (1.7,7.709371087427107) node[anchor=north west] {$A$};
\draw [line width=1.3pt] (2.5,7.5)-- (4.,8.);
\draw [line width=1.3pt] (4.,8.)-- (3.5,7.5);
\draw [line width=1.3pt] (4.,8.)-- (5.5,7.5);
\draw [line width=1.3pt] (4.5,7.5)-- (4.,8.);
\draw [line width=1.3pt] (4.,8.)-- (4.,8.5);
\draw [line width=1.3pt] (4.,8.5)-- (4.,9.);
\draw [line width=1.3pt] (4.,9.)-- (4.,9.5);
\draw (-0.106576701308053993,0.20005064773549477) node[anchor=north west] {$A$};
\draw (3.3289078315441382,6.681772291297869) node[anchor=north west] {$C_1$};
\draw (3.3289078315441382,9.169705851453124) node[anchor=north west] {$C_2$};
\begin{scriptsize}
\draw [fill=black] (2.5,7.5) circle (2.5pt);
\draw [fill=black] (3.5,7.5) circle (2.5pt);
\draw [fill=black] (4.5,7.5) circle (2.5pt);
\draw [fill=black] (5.5,7.5) circle (2.5pt);
\draw [fill=black] (4.,6.) circle (2.5pt);
\draw [fill=black] (4.,5.5) circle (2.5pt);
\draw [fill=black] (4.,6.5) circle (2.5pt);
\draw [fill=black] (4.,7.) circle (2.5pt);
\draw [fill=black] (4.,8.) circle (2.5pt);
\draw [fill=black] (4.,8.5) circle (2.5pt);
\draw [fill=black] (4.,9.) circle (2.5pt);
\draw [fill=black] (4.,9.5) circle (2.5pt);
\end{scriptsize}
\end{tikzpicture}

%% file: posetconst3.tex
\begin{tikzpicture}[line cap=round,line join=round,>=triangle 45,x=1.0cm,y=1.0cm]
\clip(1.538528068135135,3.887909507184503) rectangle (6.118394068776089,8.13962311525619);
\draw [shift={(2.5,4.5)},line width=1.3pt]  plot[domain=1.5707963267948966:3.141592653589793,variable=\t]({1.*0.3*cos(\t r)+0.*0.3*sin(\t r)},{0.*0.3*cos(\t r)+1.*0.3*sin(\t r)});
\draw [shift={(2.5,4.5)},line width=1.3pt]  plot[domain=3.141592653589793:4.71238898038469,variable=\t]({1.*0.3*cos(\t r)+0.*0.3*sin(\t r)},{0.*0.3*cos(\t r)+1.*0.3*sin(\t r)});
\draw [shift={(5.5,4.5)},line width=1.3pt]  plot[domain=-1.5707963267948966:0.,variable=\t]({1.*0.3*cos(\t r)+0.*0.3*sin(\t r)},{0.*0.3*cos(\t r)+1.*0.3*sin(\t r)});
\draw [shift={(5.5,4.5)},line width=1.3pt]  plot[domain=0.:1.5707963267948966,variable=\t]({1.*0.3*cos(\t r)+0.*0.3*sin(\t r)},{0.*0.3*cos(\t r)+1.*0.3*sin(\t r)});
\draw [line width=1.3pt] (2.2,4.5)-- (2.2,4.5);
\draw [line width=1.3pt] (2.5,4.8)-- (5.5,4.8);s
\draw [line width=1.3pt] (5.8,4.5)-- (5.8,4.5);
\draw [line width=1.3pt] (5.5,4.2)-- (2.5,4.2);
\draw [shift={(3.5,6.)},line width=1.3pt]  plot[domain=1.5707963267948966:3.141592653589793,variable=\t]({1.*0.3*cos(\t r)+0.*0.3*sin(\t r)},{0.*0.3*cos(\t r)+1.*0.3*sin(\t r)});
\draw [shift={(3.5,6.)},line width=1.3pt]  plot[domain=3.141592653589793:4.71238898038469,variable=\t]({1.*0.3*cos(\t r)+0.*0.3*sin(\t r)},{0.*0.3*cos(\t r)+1.*0.3*sin(\t r)});
\draw [shift={(5.5,6.)},line width=1.3pt]  plot[domain=-1.5707963267948966:0.,variable=\t]({1.*0.3*cos(\t r)+0.*0.3*sin(\t r)},{0.*0.3*cos(\t r)+1.*0.3*sin(\t r)});
\draw [shift={(5.5,6.)},line width=1.3pt]  plot[domain=0.:1.5707963267948966,variable=\t]({1.*0.3*cos(\t r)+0.*0.3*sin(\t r)},{0.*0.3*cos(\t r)+1.*0.3*sin(\t r)});
\draw [line width=1.3pt] (3.2,6.)-- (3.2,6.);
\draw [line width=1.3pt] (3.5,6.3)-- (5.5,6.3);
\draw [line width=1.3pt] (5.8,6.)-- (5.8,6.);
\draw [line width=1.3pt] (5.5,5.7)-- (3.5,5.7);
\draw [shift={(2.5,7.5)},line width=1.3pt]  plot[domain=1.5707963267948966:3.141592653589793,variable=\t]({1.*0.3*cos(\t r)+0.*0.3*sin(\t r)},{0.*0.3*cos(\t r)+1.*0.3*sin(\t r)});
\draw [shift={(2.5,7.5)},line width=1.3pt]  plot[domain=3.141592653589793:4.71238898038469,variable=\t]({1.*0.3*cos(\t r)+0.*0.3*sin(\t r)},{0.*0.3*cos(\t r)+1.*0.3*sin(\t r)});
\draw [shift={(5.5,7.5)},line width=1.3pt]  plot[domain=-1.5707963267948966:0.,variable=\t]({1.*0.3*cos(\t r)+0.*0.3*sin(\t r)},{0.*0.3*cos(\t r)+1.*0.3*sin(\t r)});
\draw [shift={(5.5,7.5)},line width=1.3pt]  plot[domain=0.:1.5707963267948966,variable=\t]({1.*0.3*cos(\t r)+0.*0.3*sin(\t r)},{0.*0.3*cos(\t r)+1.*0.3*sin(\t r)});
\draw [line width=1.3pt] (2.2,7.5)-- (2.2,7.5);
\draw [line width=1.3pt] (2.5,7.8)-- (5.5,7.8);
\draw [line width=1.3pt] (5.8,7.5)-- (5.8,7.5);
\draw [line width=1.3pt] (5.5,7.2)-- (2.5,7.2);
\draw [line width=1.3pt] (2.5,4.5)-- (2.5,5.);
\draw [line width=1.3pt] (2.5,5.)-- (4.5,4.5);
\draw [line width=1.3pt] (2.5,5.)-- (5.5,4.5);
\draw [line width=1.3pt] (2.5,5.)-- (3.5,4.5);
\draw [line width=1.3pt] (2.5,7.)-- (2.5,7.5);
\draw [line width=1.3pt] (2.5,7.)-- (3.5,7.5);
\draw [line width=1.3pt] (2.5,7.)-- (4.5,7.5);
\draw [line width=1.3pt] (2.5,7.)-- (5.5,7.5);
\draw (1.5382730027425286,4.786634180221476) node[anchor=north west] {$A_1$};
\draw (1.5382730027425286,7.782808199332397) node[anchor=north west] {$A_2$};
\draw (1.5382730027425286,6.282808199332397) node[anchor=north west] {$C$};
\draw (5.6882730027425286,6.282808199332397) node[anchor=north west] {$B$};
\draw [line width=1.3pt] (2.5,5.)-- (2.5,6.);
\draw [line width=1.3pt] (2.5,6.)-- (2.5,7.);
\begin{scriptsize}
\draw [fill=black] (2.5,7.5) circle (2.5pt);
\draw [fill=black] (3.5,7.5) circle (2.5pt);
\draw [fill=black] (4.5,7.5) circle (2.5pt);
\draw [fill=black] (5.5,7.5) circle (2.5pt);
\draw [fill=black] (2.5,6.) circle (2.5pt);
\draw [fill=black] (3.5,6.) circle (2.5pt);
\draw [fill=black] (4.5,6.) circle (2.5pt);
\draw [fill=black] (5.5,6.) circle (2.5pt);
\draw [fill=black] (2.5,4.5) circle (2.5pt);
\draw [fill=black] (3.5,4.5) circle (2.5pt);
\draw [fill=black] (4.5,4.5) circle (2.5pt);
\draw [fill=black] (5.5,4.5) circle (2.5pt);
\draw [fill=black] (2.5,5.) circle (2.5pt);
\draw [fill=black] (2.5,7.) circle (2.5pt);
\end{scriptsize}
\end{tikzpicture}

%% file: convex.tex
\begin{tikzpicture}[line cap=round,line join=round,>=triangle 45,x=0.35cm,y=0.35cm]
\clip(-4.449962760676213,-3.3) rectangle (11.0,10.0);
\draw [line width=2.01pt] (-3.365450203836752,2.641339866600507)-- (-3.365450203836752,6.533853071607494);
\draw [line width=2.01pt] (-3.365450203836752,6.533853071607494)-- (-1.1196718272380917,8.779631448206159);
\draw [line width=2.01pt] (-1.1196718272380917,8.779631448206159)-- (4.020951443878645,8.779631448206159);
\draw [line width=1.301pt] (4.020951443878645,8.779631448206159)-- (9.503481413642469,6.038366463324244);
\draw [line width=1.301pt] (9.503481413642469,6.038366463324244)-- (9.503481413642469,0.4906189331549795);
\draw [line width=1.301pt] (9.503481413642469,0.4906189331549795)-- (7.979855584464931,-1.03300689602256);
\draw [line width=1.301pt] (7.979855584464931,-1.03300689602256)-- (3.9832433214093754,-1.03300689602256);
\draw [line width=2.01pt] (3.9832433214093754,-1.03300689602256)-- (-3.365450203836752,2.641339866600507);
\draw [line width=2.01pt] (5.0,5.0)-- (4.020951443878645,8.779631448206159);
\draw [line width=2.01pt] (5.0,5.0)-- (3.9832433214093754,-1.03300689602256);
\draw [line width=2.0pt,dash pattern=on 3pt off 5pt] (3.9832433214093754,-1.03300689602256)-- (4.020951443878645,8.779631448206159);

\draw (5.505318994531624,5.0162293412155545) node[anchor=north west] {$P$};
\draw (2.2,-0.8) node[anchor=north west] {$v_0$};
\draw (7.020840963777292,0.5647014967864794) node[anchor=north west] {$v_7$};
\draw (8.0,1.9627057007042479) node[anchor=north west] {$v_6$};
\draw (9.6,7.3) node[anchor=north west] {$v_5$};
\draw (3.8,10.1) node[anchor=north west] {$v_4$};
\draw (-2.3,10.1) node[anchor=north west] {$v_3$};
\draw (-4.6,7.94903139440341) node[anchor=north west] {$v_2$};
\draw (-4.6,2.6437846718436733) node[anchor=north west] {$v_1$};
\begin{scriptsize}
\draw [fill=black] (3.9832433214093754,-1.03300689602256) circle (2.5pt);
\draw [fill=black] (7.979855584464931,-1.03300689602256) circle (2.5pt);
\draw [fill=black] (9.503481413642469,0.4906189331549795) circle (2.5pt);
\draw [fill=black] (9.503481413642469,6.038366463324244) circle (2.5pt);
\draw [fill=black] (4.020951443878645,8.779631448206159) circle (2.5pt);
\draw [fill=black] (-1.1196718272380917,8.779631448206159) circle (2.5pt);
\draw [fill=black] (-3.365450203836752,6.533853071607494) circle (2.5pt);
\draw [fill=black] (-3.365450203836752,2.641339866600507) circle (2.0pt);
\draw [fill=black] (5.0,5.0) circle (2.5pt);
\end{scriptsize}
\end{tikzpicture}\qquad
\begin{tikzpicture}[line cap=round,line join=round,>=triangle 45,x=0.35cm,y=0.35cm]
\clip(-4.449962760676213,-3.3) rectangle (14.0,10.0);
\draw [line width=2.01pt] (-3.365450203836752,2.641339866600507)-- (-3.365450203836752,6.533853071607494);
\draw [line width=2.01pt] (-3.365450203836752,6.533853071607494)-- (-1.1196718272380917,8.779631448206159);
\draw [line width=2.01pt] (-1.1196718272380917,8.779631448206159)-- (4.020951443878645,8.779631448206159);
\draw [line width=2.01pt] (4.020951443878645,8.779631448206159)-- (9.503481413642469,6.038366463324244);
\draw [line width=1.301pt] (9.503481413642469,6.038366463324244)-- (9.503481413642469,0.4906189331549795);
\draw [line width=1.301pt] (9.503481413642469,0.4906189331549795)-- (7.979855584464931,-1.03300689602256);
\draw [line width=1.301pt] (7.979855584464931,-1.03300689602256)-- (3.9832433214093754,-1.03300689602256);
\draw [line width=2.01pt] (3.9832433214093754,-1.03300689602256)-- (-3.365450203836752,2.641339866600507);
\draw [line width=2.01pt] (12.15471204908579,-2.4388474239735722)-- (3.9832433214093754,-1.03300689602256);
\draw [line width=2.01pt] (12.15471204908579,-2.4388474239735722)-- (9.503481413642469,6.038366463324244);
\draw (12.505318994531624,-2.0162293412155545) node[anchor=north west] {$P$};
\draw (2.2,-0.8) node[anchor=north west] {$v_0$};
\draw (7.020840963777292,0.5647014967864794) node[anchor=north west] {$v_7$};
\draw (8.0,1.9627057007042479) node[anchor=north west] {$v_6$};
\draw (9.6,7.3) node[anchor=north west] {$v_5$};
\draw (3.8,10.1) node[anchor=north west] {$v_4$};
\draw (-2.3,10.1) node[anchor=north west] {$v_3$};
\draw (-4.6,7.94903139440341) node[anchor=north west] {$v_2$};
\draw (-4.6,2.6437846718436733) node[anchor=north west] {$v_1$};
\begin{scriptsize}
\draw [fill=black] (3.9832433214093754,-1.03300689602256) circle (2.5pt);
\draw [fill=black] (7.979855584464931,-1.03300689602256) circle (2.5pt);
\draw [fill=black] (9.503481413642469,0.4906189331549795) circle (2.5pt);
\draw [fill=black] (9.503481413642469,6.038366463324244) circle (2.5pt);
\draw [fill=black] (4.020951443878645,8.779631448206159) circle (2.5pt);
\draw [fill=black] (-1.1196718272380917,8.779631448206159) circle (2.5pt);
\draw [fill=black] (-3.365450203836752,6.533853071607494) circle (2.5pt);
\draw [fill=black] (-3.365450203836752,2.641339866600507) circle (2.0pt);
\draw [fill=black] (12.15471204908579,-2.4388474239735722) circle (2.5pt);
\end{scriptsize}
\end{tikzpicture}